\documentclass[10pt]{article}

\usepackage[ruled,vlined,linesnumbered]{algorithm2e}
\usepackage{subcaption}
\usepackage{amsmath}
\usepackage{mathtools}
\usepackage{amssymb}
\usepackage{amsthm}
\usepackage{thmtools}
\usepackage{thm-restate}
\usepackage{comment}
\usepackage{color} 
\usepackage[unicode]{hyperref}
\usepackage[T1]{fontenc}
\usepackage{geometry}
\usepackage{multirow}
\usepackage[colorinlistoftodos]{todonotes}
\usepackage{lmodern}
\usepackage{anyfontsize}
\usepackage{stmaryrd}
\usepackage{natbib}
\usepackage{cleveref}
 \usepackage{graphicx}
\usepackage{subcaption}
\usepackage{dsfont}
\usepackage{float} 
\usepackage{amsbsy}
\usepackage{sidecap}
\usepackage{cancel}
\usepackage{caption}
\usepackage{subcaption}
\setcitestyle{numbers,open={[},close={]}}
\definecolor{red}{rgb}{0.7,0.15,0.15}
\definecolor{green}{rgb}{0,0.5,0}
\definecolor{blue}{rgb}{0,0,0.7}
\hypersetup{colorlinks, linkcolor={red},citecolor={green}, urlcolor={blue}}
			
\makeatletter \@addtoreset{equation}{section}

\newtheorem{theorem}{Theorem}[section]

\newtheorem{corollary}[theorem]{Corollary}

\newtheorem{lemma}[theorem]{Lemma}
\newtheorem{proposition}[theorem]{Proposition}

\newtheorem{definition}[theorem]{Definition}
\newtheorem{remark}[theorem]{Remark}

\def\no{\noindent}
\def\beq{\begin{eqnarray}}
\def\eeq{\end{eqnarray}}
\def\be*{\begin{eqnarray*}}
\def\ee*{\end{eqnarray*}}

\DeclareUnicodeCharacter{014D}{\=o}
 \title{Existence of a non-zero steady-state solution to a population density equation subject to horizontal transfer}

\author{Mateo Deangeli Bravo\thanks{Ecole Polytechnique Paris, Centre de Math\'ematiques Appliqu\'ees,
        mateo.deangeli@polytechnique.edu.}
    }
             \date{\today}

\begin{document}

\maketitle
 
\begin{abstract}
From plasmid dynamics and efflux pump expression in bacteria to risk-taking behaviour in human groups, horizontal interactions can drive the transmission of deleterious traits that spread despite carrying no selective advantage, thereby burdening the population. How can such detrimental traits persist, and under what conditions can the population avoid extinction?

To address these questions, we model a population structured by a continuous trait representing individual burden within a non-linear, non-local integro-differential framework. From an analytical perspective, mutations are modeled via a jump operator rather than a diffusion term; consequently, the mutation operator lacks regularizing properties and is generally non-self-adjoint. Moreover, working in a weighted $L^1$ setting precludes the direct application of standard Krein–Rutman theory. To overcome these hurdles, we construct a custom fixed-point operator by first freezing the non-linearities to solve an auxiliary eigenvalue problem, and then tuning the resulting principal eigenvalue. Applying the Leray–Schauder fixed-point theorem to this operator rigorously yields the existence of a non-trivial steady state. Finally, numerical simulations illustrate the qualitative dynamics of the model and delineate the parameter regimes governing population persistence and extinction.
\end{abstract}

\vspace{3mm}
\no{\bf Keywords.} Integros-differential equations, Non-local jump operators, Leray–Schauder fixed-point theorem, Numerical analysis of local stability, Horizontal gene transfer, Population persistence. 

\vspace{3mm}
\no{\bf MSC2010.} 35R09, 47H10, 47J10, 92D15.

\section{Introduction}
Horizontal gene transfer (HGT) is a fundamental driver of bacterial diversification \cite{dewan2023mathematician}, notably through plasmid conjugation. While plasmids can carry adaptive genes, expressing them imposes a metabolic cost. This mechanism generalizes naturally to continuous biological traits, such as the expression of efflux pumps. Through quorum sensing, cells with higher pump levels trigger a positive feedback loop in neighbouring bacteria, driving them toward higher expression thresholds that carry heavy metabolic costs in neutral environments. Beyond microbiology, this asymmetric alignment also applies to social dynamics and the spread of harmful behaviours, such as smoking intensity or exposure to risky behaviours. Here, "transfer" occurs through peer pressure or social imitation, pushing individuals toward higher risk levels without providing any evolutionary benefit. To translate these mechanisms into mathematical terms, we consider that transfer depends on the relative load between individuals. Encounters are intrinsically directional.

To study these dynamics, we model the evolution of a population of individuals, each characterized by a real number $x \in \mathbb{R}_+$ representing the intensity of the burden (e.g., the quantity of plasmid DNA, efflux pump expression level, or the intensity of a risky behavioural trait). We make the hypothesis that the environment is neutral, meaning that the trait confers no selective advantage and represents a pure maintenance or exposure cost. The population density $n(t,x)$ of individuals with trait $x$ at time $t$ is given by the solution of an integro-differential equation—the existence of which can be rigorously obtained as the limit of a random individual-based processes (see \cite{bansaye2015stochastic, A2, deangelibravo2026ancestral, fournier_meleard}):
\begin{equation}\label{eq:integrodiff}
    \begin{aligned} 
\forall (t,x) \in \mathbb{R}^2_+, \quad & \\ 
\partial_t n(t,x) =& \left(r(x)-C\int_0^\infty n(t,y)dy\right)n(t,x)+\int_0^\infty \left[\gamma(y)n(t,y)m(y,x)-\gamma(x)n(t,x)m(x,y)\right]dy\\ 
&+\frac{\tau n(t,x)}{\beta+\mu \int_0^\infty n(t,y)dy}\int_0^\infty \text{sgn}(x-y)n(t,y)dy
\end{aligned}
\end{equation}

where $r(x)$ is the growth term of an individual with trait $x$, $C$ the competition intensity, here assumed constant, $\gamma(x)$ the frequency of occurrence of trait variation events (e.g., mutation, trait loss, or behavioural shift) for an individual with trait $x$, $m$ the variation kernel that can be heavy-tailed. For the transfer, the individual with the larger trait drives the exchange toward the recipient. We capture this asymmetric interaction using the sign function, $\text{sgn}(x-y)$, which sets the direction of the transfer flow between individuals with traits $x$ and $y$. The parameters for the transfer term are $\tau$, the transfer intensity, $\beta$ the encounter time between two individuals, and $\mu$ the average time taken by the transfer process. Depending on the values of these parameters, the transfer regime varies. It is described as Frequency-Dependent when $\beta=0$ and as Density-Dependent when $\mu=0$. When both $\mu$ and $\beta$ are strictly positive, it is referred to as a mixed regime, called the Beddington-DeAngelis regime (see \cite{A1,A2}).
\begin{figure}[H]
    \centering
    \begin{subfigure}[b]{0.35\textwidth}
        \centering
        \includegraphics[width=\linewidth]{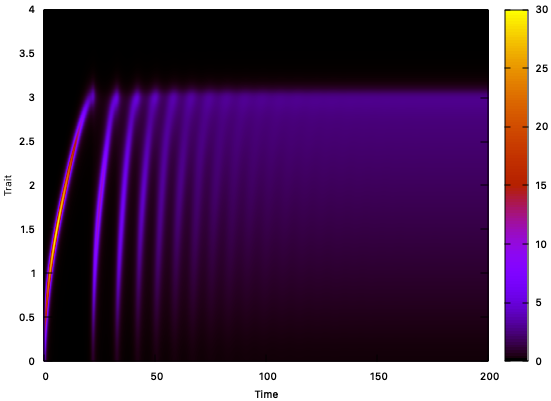}
        \caption{\small Convergence to a stationary state - Beddington DeAngelis transfer case}
        \label{fig:image1}
    \end{subfigure}
    \hspace{0.05\textwidth}
    \begin{subfigure}[b]{0.35\textwidth}
        \centering
        \includegraphics[width=\linewidth]{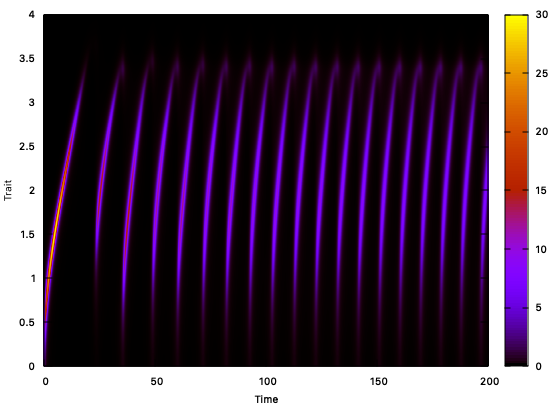}
        \caption{\small Slow convergence oscillations - Frequency Dependence transfer case}
        \label{fig:image2}
    \end{subfigure}

    \caption{Numerical analysis of Eq. \eqref{eq:integrodiff}: Time evolution of the density $n$ showing (a) convergence in the BDA regime and (b) undamped oscillations in the FD regime $(\beta=\mu=1)$.}
    \label{fig:dynamics}
\end{figure}

Solving the equation numerically reveals convergence to a positive stationary profile in non-frequency-dependent regimes (such as density-dependent or Beddington-DeAngelis). Conversely, in the frequency-dependent case, moderately large values of $\tau$ induce long-lived transient oscillations that decay very slowly toward a positive equilibrium. It brings some light on individual-based simulations of \cite{A2}.
In contrast, if $\tau$ is too large, the system bypasses these prolonged fluctuations and directly collapses to extinction, a behaviour that arises precisely when the existence of a non-zero stationary solution is not theoretically guaranteed. In related models, notably with different rescaling, rigorous proofs of cyclic behavior have only been obtained in \cite{champagnat2021stochastic}, while it is conjectured in \cite{A3}. Furthermore, establishing the stability of the stationary regimes remains an open challenge.

As a first step in the study of this equation, we focus here on the existence of a non-trivial stationary solution $\overline{n}$ :
\begin{equation}
    \begin{aligned} 
&(r(x)-C\rho)\overline{n}(x)+ \int_0^\infty[\gamma(y) \overline{n}(y)m(y,x)-\gamma(x)\overline{n}(x)m(x,y)]dy+\overline{n}(x)\overline{\mathcal{H}}(\overline{n},\rho)(x)=0, \\ 
&\rho=\int_0^\infty \overline{n}(y)dy > 0, \\ 
&\overline{\mathcal{H}}(\overline{n},\rho)(x)=\frac{\tau}{\beta+\mu\rho}\int_0^\infty \text{sgn}(x-y)\overline{n}(y)dy, \\ 
&\overline{n}(x) \geq 0, \quad \int_0^\infty \overline{n}(y)dy>0.
\end{aligned}
\end{equation}

Assuming the solution is not trivial, it can then be renormalised to obtain a probability density $u$. We thus seek a pair $(u,\rho)$ such that:

\begin{equation}\label{eq:stationnaire}
    \begin{aligned} 
&(r(x)-C\rho)u(x)+ \int_0^\infty[\gamma(y) u(y)m(y,x)-\gamma(x)u(x)m(x,y)]dy+u(x)\mathcal{H}(u,\rho)(x)=0, \\ 
&\rho=\frac{1}{C}\int_0^\infty r(y)u(y)dy > 0, \\ 
&\mathcal{H}(u,\rho)(x)=\frac{\tau \rho}{\beta+\mu\rho}\int_0^\infty \text{sgn}(x-y)u(y)dy, \\ 
&u(x) \geq 0, \quad \int_0^\infty u(y)dy=1. 
\end{aligned}
\end{equation}
From this pair $(u,\rho)$, we can reconstruct $\overline{n}(x)=\rho u(x)$. The parameter $\rho$ is defined by integrating the first equation of \eqref{eq:stationnaire}. In the unrenormalised case, it represents the population size, whereas in the renormalised equation, it corresponds to the overall fitness of the population. 
Notice that the transfer operator $\mathcal{H}(u,\rho)$ involves the non-linear intensity factor $\frac{\tau \rho}{\beta+\mu\rho}$, which encapsulates two distinct biological limits. When the population is near extinction ($\rho \to 0$), the transfer is limited by the time needed to find a partner ($\beta$), yielding a linear behavior $\frac{\tau}{\beta}\rho$. Conversely, for large $\rho$, the transfer saturates at $\frac{\tau}{\mu}$, limited by the duration of the transfer process itself ($\mu$). This distinction becomes mathematically crucial in the pure frequency-dependent regime ($\beta=0$): in this case, partner encounters are instantaneous, and the transfer intensity remains constantly equal to $\tau/\mu$ regardless of population density. Because this deleterious effect persists even at infinitesimally small densities, it imposes a strictly stronger condition for the existence of a stationary state (typically requiring a spectral bound $\lambda_0 > \tau/\mu$ rather than simply $\lambda_0 > 0$).
\begin{figure}[H]
    \centering
    \begin{subfigure}[b]{0.35\textwidth}
        \centering
        \includegraphics[width=\linewidth]{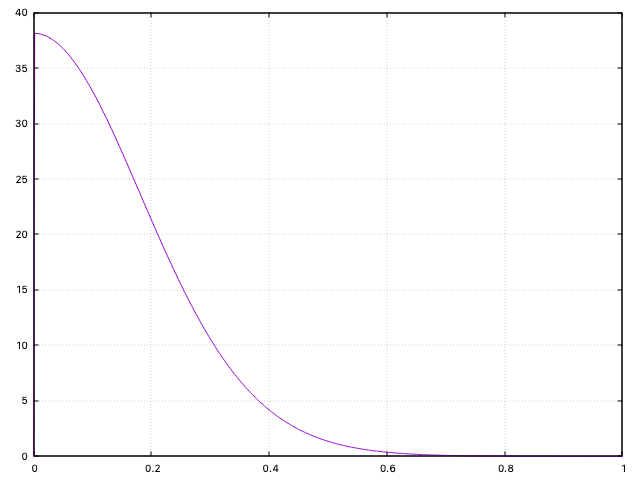}
        \caption{\small The stationary solution ($\tau=0$)}
        \label{fig:stationary_without_transfer}
    \end{subfigure}
    \hspace{0.05\textwidth}
    \begin{subfigure}[b]{0.35\textwidth}
        \centering
        \includegraphics[width=\linewidth]{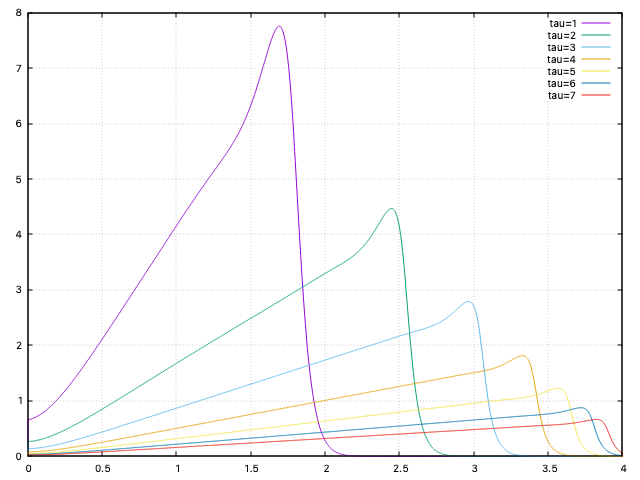}
        \caption{\small Possible stationary solution for $\tau=1,...,7$}
        \label{fig:stationaries_with_transfer}
    \end{subfigure}
    \caption{Steady-state profiles (a) without transfer and (b) for different values of $\tau$ in the BDA regime.}
    \label{fig:profiles}
\end{figure}
To prove the existence of this pair $(u, \rho)$, we build upon the approach initiated in \cite{mirrahimi} using the Leray-Schauder fixed-point theorem \cite{mawhin1999leray}. In \cite{mirrahimi}, the authors investigated a related model driven by horizontal gene transfer, with a main focus on deriving its Hamilton-Jacobi limit. While we share a similar biological motivation, our objective and mathematical setting differ, giving rise to several key technical difficulties. First, mutations and variations in trait value are modeled here by a jump operator rather than a classical diffusion term. As a direct consequence, the mutation operator completely lacks parabolic regularizing properties. Furthermore, our analysis operates in a weighted $L^1$ setting, which rules out the standard compactness arguments required to apply the classical Krein–Rutman theorem. Compounding this issue, our framework entirely bypasses the self-adjointness assumption on the linearized operator without transfer—a property that played a pivotal role in the spectral analysis of \cite{mirrahimi}. In addition, we consider a more general horizontal transfer intensity operator capable of encompassing the broad spectrum of biological and social phenomena mentioned above, thereby introducing stronger and less tractable non-linearities. To overcome these intertwined analytical hurdles, we design a fixed-point scheme. Instead of inverting the system, we first freeze the non-linearities to formulate an auxiliary, linearized eigenvalue problem in our weighted space. Once the existence of a principal eigenvalue and an associated positive eigenfunction is established, we calibrate this eigenvalue so as to recover the correct non-linear scaling. The fixed-point mapping defined through this two-step procedure is then shown to satisfy the compactness and topological requirements necessary to apply the Leray–Schauder theorem, ultimately yielding the existence of a non-trivial steady state.

We obtain the existence of a solution to equation \eqref{eq:stationnaire} by extending the existence result from the model without transfer, namely:

\begin{equation}\label{eq:stationnaire_sans_transfert}
    \begin{aligned} 
&(r(x)-C\rho_0)N(x)+ \int_0^\infty[\gamma(y) N(y)m(y,x)-\gamma(x)N(x)m(x,y)]dy=0, \\ 
&\rho_0=\frac{1}{C}\int_0^\infty r(y)N(y)dy > 0, \\ 
&N(x) \geq 0, \quad \int_0^\infty N(y)dy=1. 
\end{aligned}
\end{equation}

Our main result reads as follows.

\begin{theorem}
Equation \eqref{eq:stationnaire_sans_transfert} admits a solution. Assume moreover that
\begin{equation}\label{eq:cond_rho0}
C\rho_0 > \frac{\tau}{\mu} \mathds{1}_{\{\beta=0\}}.
\end{equation}
Then, there exists a solution to \eqref{eq:stationnaire}.
\end{theorem}
The assumptions that ensure the existence and uniqueness to \eqref{eq:stationnaire_sans_transfert} are stated below and the proof is given in Section \ref{section:pb_without_transfer}.

Numerically, the sufficient condition \eqref{eq:cond_rho0} for a non-trivial stationary profile fails precisely when the integro-differential equation converges to the zero state. While the present work rigorously settles the question of existence, establishing the uniqueness and global asymptotic stability of these non-trivial solutions remains an open challenge. Despite our efforts to find an analytical path, standard tools such as General Relative Entropy methods (see \cite{michel2007general, michel2005general, perthame2007transport}) ultimately fall short: the different natures of the linear variation operator and the non-local nonlinear transfer term prevent the derivation of a common Lyapunov functional. To overcome this theoretical obstacle, we turn to a numerical study to assess the linear stability of the profiles and capture their long-time convergence dynamics.

We begin by specifying the assumptions with which we are working. Section \ref{section:pb_without_transfer} proves the existence of a solution to the steady-state equation without the transfer term, necessary for the application of the Leray-Schauder theorem. The resolution of this equation also allows one to be familiar with the definition of the fixed-point operator presented in Section \ref{section:pb_with_transfer}. We then demonstrate the existence of a steady-state solution to the problem with transfer under a condition related to the eigenvalue of the problem without transfer. Finally, in Section \ref{section:explicit_conditions} we present explicit conditions for existence and how the parameters of the model interact with one another.
Continuity and compactness of the involved operators are postponed in \ref{section:regularity_results} to improve readibility.
Lastly, we study numerically the convergence of the integrodifferential equation and the stability of the steady-states in Section \ref{sec:numerical_study}.

\section*{Assumption \textbf{(H)}}
\begin{itemize}
    \item (H.a) $r$ is continuous, decreasing, strictly positive at $0$, and $\lim_{x\rightarrow \infty}r(x)=-\infty$.
    \item (H.b) The kernel $m$ is strictly positive, uniformly continuous in both variables and
    \begin{equation*}
        \int_0^\infty m(x,y)dy=1,\qquad \overline{M}:=\sup_{y\geq0}\int_0^\infty m(x,y)dx<\infty.
    \end{equation*}
    \item (H.c) The mutation rate $\gamma$ is continuous, bounded by $\overline{\gamma}$, strictly positive almost everywhere and 
    \begin{equation}\label{eq:condition_forte_de_positivite}
         \sup_{x\geq0}(r(x)-\gamma(x))>0.
    \end{equation}
    \item (H.d) $C$ and $\tau$ are positive, $\beta$ and $\mu$ are non negative and non-zero simultaneously.
\end{itemize}
The fact that $r$ diverges to negative infinity ensures that the stationary solution is bounded and we need $r$ to be continuous to obtain the continuity of the solution. The strict positivity of the kernel $m$ guarantees that the principal eigenvalue is strictly positive, while its uniform continuity ensures the continuity of the operators under consideration. The assumption \eqref{eq:condition_forte_de_positivite} ensures that the unique eigenvalue for the problem without transfer is strictly positive. It is not absolutely necessary but is convenient when $\gamma$ is not constant and $m$ is not symmetric. Biologically, it can be justified by the idea that the maximum reproduction rate is greater than the variation rate, which is particularly clear when $\gamma(x)=\overline{\gamma}$ since it can be written as $r_0=r(0)>\overline{\gamma}$.

\subsection*{Notations :}
We mainly work in the space $L^1(d\nu)$ defined by:
\begin{equation*}
\begin{aligned}
    &\nu(dx):=(1+|r(x)|)dx,\\
    &L^1(d\nu):= \{u:\mathbb{R}_+\rightarrow\mathbb{R}\,|\, \int_0^\infty |u(x)|\nu(dx)<\infty\},\\
    &\|u\|_\nu:=\int_0^\infty |u(x)|\nu(dx).
\end{aligned}
\end{equation*}
We write $L^1(dx)$ for the space of integrable functions with respect to the Lebesgue measure on $\mathbb{R}_+$.
\section{The Transfer-Free Equation}\label{section:pb_without_transfer}

We first prove the existence and uniqueness of the solution to the equation without the transfer term \eqref{eq:stationnaire_sans_transfert}. Then, we demonstrate the existence and uniqueness of the dual problem (without the transfer term), which will be used for establishing a priori estimates.

We define
\begin{equation}\label{eq:def_linear_operator}
    \begin{aligned}
        L\phi(x)&=r(x)\phi(x)+\int_0^\infty (\gamma(y)\phi(y)m(y,x)-\gamma(x)\phi(x)m(x,y))dy\\
        &=(r(x)-\gamma(x))\phi(x)+\int_0^\infty \gamma(y)\phi(y)m(y,x)dy
    \end{aligned}
\end{equation}
where the form given in the second line follows from (H.b).
\begin{proposition}\label{prop:existence_sol_sans_transfert}
Under the assumption \textbf{(H)}, there exists a unique pair $(N,\lambda_0)\in L^1(d\nu)\times\mathbb{R}^*_+$ such that
\begin{equation}\label{eq:pb_sans_transfert}
    \begin{aligned}
    &LN(x)=\lambda_0 N(x),\, x\in\mathbb{R}_+\\
    &N(x)>0,\, x\in\mathbb{R}_+\\
    &\int_0^\infty N(x)dx=1.
\end{aligned}
\end{equation}
Moreover, $N$ is a continuously bounded function.
\end{proposition}
\begin{remark}
The equation \eqref{eq:pb_sans_transfert} is equivalent to \eqref{eq:stationnaire_sans_transfert} with $\lambda_0=C\rho_0$.
\end{remark}
The proof is based on rewriting problem \eqref{eq:pb_sans_transfert} as a parameterized fixed point problem. We first show the existence of a solution for the whole range of parameters $\lambda$ using a Krein-Rutman-type Theorem, and we conclude by proving that there is a unique parameter corresponding to \eqref{eq:pb_sans_transfert}.
\begin{proof}
We rewrite equation \eqref{eq:pb_sans_transfert} in the form of a fixed-point problem. We define
\begin{equation*}
    D_\lambda(x):= \lambda+\gamma(x)-r(x)
\end{equation*}
and 
\begin{equation*}
    \Sigma:=\{\lambda\in\mathbb{R}|\, \forall x\geq0,\,D_\lambda(x)>0\}=\left\{\lambda\in\mathbb{R}\Big|\lambda> \sup_{x\geq0}(r(x)-\gamma(x))\right\}
\end{equation*}
the set of admissible values for the parameter $\lambda$. By Assumption (H.c), $\Sigma\subset\mathbb{R}^*_+$. For any $\lambda\in\Sigma$ and any $\phi\in L^1(d\nu)$, we define the operator $T_\lambda$ evaluated at $x \geq 0$ by
\begin{equation*}
    T_\lambda \phi(x):=\frac{1}{D_\lambda(x)}\int_0^\infty\gamma(y)m(y,x)\phi(y)dy.
\end{equation*}
The equation \eqref{eq:pb_sans_transfert} is equivalent to $T_{\lambda_0}N=N$.

The value of $\lambda_0$ is unknown and will be determined later. We will show that for any $\lambda \in \Sigma$, there exists a pair $(\phi_\lambda, \psi(\lambda))\in L^1(d\nu)\times\mathbb{R}_+^*$ such that 
\begin{equation}\label{eq:auxiliary_operator_eigenvalue}
T_\lambda \phi_\lambda = \psi(\lambda)\phi_\lambda.
\end{equation}
It will then suffice to find the unique $\lambda \in \Sigma$ satisfying $\psi(\lambda)=1$ to recover $T_{\lambda_0}N=N$.

\underline{Existence of the pair $(\phi_\lambda, \psi(\lambda))$:}

We use Jentzsch's Theorem (Th.\ref{thm:Jentzsch} in Appendix \ref{appendix}) to obtain \eqref{eq:auxiliary_operator_eigenvalue}. The proof of the compactness of the operator $T_\lambda$ on $L^1(d\nu)$ for any $\lambda$ in $\Sigma$ is postponed to Lemma \ref{lemma:compactness_of_operator_without_transfer}. To apply Theorem \ref{thm:Jentzsch}, we consider the $\sigma$-finite measured space $(\mathbb{R}_+, \mathcal{B}(\mathbb{R}_+), \nu)$. The operator $T_\lambda : L^1(d\nu) \to L^1(d\nu)$ is linear and compact, and its integral action with respect to the measure $\nu$ reads:
\begin{equation*}
    T_\lambda \phi(x) = \int_0^\infty \tilde{K}(x,y) \phi(y) d\nu(y),
\end{equation*}
where the kernel is defined by
\begin{equation*}
    \tilde{K}(x,y):=\frac{\gamma(y)m(y,x)}{D_\lambda(x)(1+|r(y)|)}.
\end{equation*}
By Assumptions (H.b) and (H.c), $\tilde{K}$ is almost everywhere positive. For all $A \in \mathcal{B}(\mathbb{R}_+)$ such that $\nu(A) > 0$ and $\nu(\mathbb{R}_+\backslash A) > 0$, we have
\begin{equation*}
    \int_{\mathbb{R}_+\backslash A}\int_A \tilde{K}(x,y) d\nu(x)d\nu(y)>0.
\end{equation*}
According to Theorem \ref{thm:Jentzsch}, there exists a simple eigenvalue $\psi(\lambda)$ with a unique renormalised eigenvector $\phi_\lambda \in L^1(d\nu)$ that is almost everywhere positive, and every other eigenvalue of $T_\lambda$ has a modulus strictly less than $\psi(\lambda)$. 

We now need to determine $\lambda$ such that $\psi(\lambda)=1$.

\underline{Continuity of $\psi$:}

To establish the continuity of $\psi$, we verify the assumptions of Theorem \ref{thm:continuite_du_rayon_spectral} in Appendix \ref{appendix}.

By Assumptions (H.a) and (H.c), $\Sigma$ is an open interval. Let $\lambda \in \Sigma$ and let $(\lambda_n)_{n \geq 0}$ be a sequence in $\Sigma$ converging to $\lambda$. Therefore $\Lambda := \min(\inf_n \lambda_n, \lambda)$ is an element of $\Sigma$. Given $f \in L^1(d\nu)$, straightforward computations yield:
$$
\begin{aligned}
    \|T_{\lambda_n} f - T_\lambda f\|_\nu &\leq \overline{\gamma} \int_0^\infty \int_0^\infty m(y,x) |f(y)| \left|\frac{1}{D_{\lambda_n}(x)} - \frac{1}{D_\lambda(x)}\right|(1+|r(x)|) dx dy \\
    &\leq \overline{\gamma} |\lambda_n - \lambda| \int_0^\infty |f(y)| \left(\int_0^\infty m(y,x) \frac{1+|r(x)|}{D_{\lambda_n}(x)D_\lambda(x)} dx\right) dy \\
    &\leq \overline{\gamma} |\lambda_n - \lambda| \int_0^\infty |f(y)| \left(\int_0^\infty m(y,x) \frac{1+|r(x)|}{D_\Lambda^2(x)} dx\right) dy.
\end{aligned}
$$
By (H.a), there exists a constant $C' > 0$ such that for all $x\geq0$, $\frac{1+|r(x)|}{D_\Lambda^2(x)} \leq C'$. Combined with Assumption (H.b), there exists a constant $C'' > 0$ such that
$$
\sup_{\|f\|_\nu \leq 1} \|T_{\lambda_n} f - T_\lambda f\|_\nu \leq C'' |\lambda_n - \lambda|.
$$
By virtue of Theorem \ref{thm:continuite_du_rayon_spectral}, the function $\lambda \mapsto \psi(\lambda)$ is continuous.

\underline{Strict decrease of $\psi$:}

We now prove the strict decrease of $\psi$ in order to ensure the uniqueness of the parameter $\lambda$ satisfying $\psi(\lambda)=1$. We apply Marek's Min-Max formula (Theorem \ref{thm:Marek}) to the operator $T_\lambda$ and first verify its assumptions. We are working in the Banach space $L^1(d\nu)$ equipped with the positive cone $L^1_+(d\nu)$, which is closed, normal, and reproducing. The operator $T_\lambda$ is compact and positive. Since the kernel $m(x,y)$ is strictly positive, the operator is \textit{semi-nonsupporting} for the integer $p=1$ (see Definition \ref{def:smnnsp}). 

We have already established via Theorem \ref{thm:Jentzsch} that the spectral radius $\psi(\lambda)$ is a simple eigenvalue admitting a positive eigenvector $\phi_\lambda \in L^1_+(d\nu)$. Furthermore, Lemma \ref{lemma:marek} ensures the existence of a positive dual eigenvector $\phi'_\lambda$ in the dual space of $L^1(d\nu)$. Consequently, $T_\lambda$ satisfies all the assumptions of Theorem \ref{thm:Marek}, and its spectral radius is characterised by:
\begin{equation}\label{eq:min_max_sup_inf}
    \psi(\lambda) = \sup_{\phi \in L^1_+(d\nu), \phi \neq 0} \left( \operatorname*{ess\,inf}_{x \in \mathbb{R}_+} \frac{T_\lambda \phi(x)}{\phi(x)} \right).
\end{equation}

Therefore, we can demonstrate the strict decrease of $\psi$ on $\Sigma$. Let $\lambda_1 < \lambda_2$ be elements of $\Sigma$. Given the strict increase of $D_\lambda$ with respect to $\lambda$, for any non-zero positive function $f \in L^1(d\nu)$, we have
\begin{equation*}
    T_{\lambda_1}f(x) > T_{\lambda_2}f(x), \quad \forall x \geq 0.
\end{equation*}
We know that there exists a unique pair $(\phi_2, \psi(\lambda_2))$ with $\phi_2 > 0$ such that
\begin{equation*}
    T_{\lambda_2}\phi_2 = \psi(\lambda_2)\phi_2,
\end{equation*}
which implies
\begin{equation*}
    T_{\lambda_1}\phi_2 > \psi(\lambda_2)\phi_2.
\end{equation*}
Applying formula \eqref{eq:min_max_sup_inf} with the operator $T_{\lambda_1}$ and the test function $\phi_2$, we obtain:
\begin{equation*}
    \psi(\lambda_1) \geq \operatorname*{ess\,inf}_{x \in \mathbb{R}_+} \frac{T_{\lambda_1}\phi_2(x)}{\phi_2(x)} > \psi(\lambda_2).
\end{equation*}
Because the underlying functions are continuous and belong to the weighted space $L^1(d\nu)$, their decay at infinity ensures that this infimum is actually attained as a minimum on $\mathbb{R}_+$. This guarantees that the strict inequality is preserved, thereby establishing the strict decrease of the mapping $\lambda \mapsto \psi(\lambda)$. 

\underline{Conclusion:}
Regarding the asymptotic behavior, as $\lambda$ tends to $+\infty$, $D_\lambda(x)$ diverges uniformly to $+\infty$, which implies that $\lim_{\lambda \to +\infty} \psi(\lambda) = 0$. Conversely, as $\lambda$ tends to $\inf\Sigma$, we have $\inf_{x \ge 0} D_\lambda(x) \to 0$. Since this minimum approaches zero, the norm of the multiplier $1/D_\lambda$ blows up, yielding $\lim_{\lambda \to \inf\Sigma} \psi(\lambda) = +\infty$. 

Since $\psi$ is continuous, the Intermediate Value Theorem guarantees the existence of a parameter $\lambda_0$ such that $\psi(\lambda_0)= 1$. Furthermore, the strict monotonicity of $\psi$ ensures that this $\lambda_0$ is unique. Consequently, the operator $T_{\lambda_0}$ admits $1$ as its principal eigenvalue. Let $N = \phi_{\lambda_0}$ denote the associated eigenfunction. Ultimately, this strict decrease also guarantees the uniqueness of the solution to equation \eqref{eq:pb_sans_transfert}, and ensures that $\lambda_0$ is both an isolated and dominant eigenvalue.

The eigenfunction $N$ is a solution of
\begin{equation*}
    N(x)=\frac{1}{\lambda_0+\gamma(x)-r(x)}\int_0^\infty \gamma(y)N(y)m(y,x)dy.
\end{equation*}
The continuity of $r$ and $\gamma$, combined with the properties of the kernel $m$, imply that $N$ is continuous. Furthermore, the fact that $-r(x)$ tends to infinity as $x \to \infty$ ensures that $N$ is bounded.

Finally, the solution of the stationary equation \eqref{eq:stationnaire} without transfer is recovered. Indeed, integrating \eqref{eq:pb_sans_transfert} over $\mathbb{R}_+$ with respect to the Lebesgue measure yields $\lambda_0 = \int_0^\infty r(y)N(y)dy$, thus reconstructing $\lambda_0 = C\rho_0$. Hypothesis (H.c) naturally ensures that $\lambda_0>0$.
\end{proof}
To establish the a priori estimates in \eqref{eq:debut_priori} needed to apply the Leray-Schauder theorem, we require a solution to the dual problem to exist.

\begin{corollary}\label{cor:existence_sans_transfert_duale}
Let $L^*$ denote the formal adjoint operator of $L$ with respect to the Lebesgue measure on $\mathbb{R}_+$. For all $\phi_1,\phi_2 \in L^1(d\nu) \cap C_b(\mathbb{R}_+)$, we have
\begin{equation*}
    \int_0^\infty L\phi_1(x) \phi_2(x)dx = \int_0^\infty \phi_1(x) L^*\phi_2(x)dx,
\end{equation*}
where
\begin{equation*}
    L^*\phi_2(x)=r(x)\phi_2(x)+\gamma(x)\int_0^\infty (\phi_2(y)-\phi_2(x))m(x,y)dy.
\end{equation*}

Under the same assumptions as Proposition \ref{prop:existence_sol_sans_transfert}, there exists a unique dual eigenvector $\varphi \in L^1(d\nu)\cap C_b(\mathbb{R}_+)$ such that 
\begin{equation}\label{eq:pb_sans_transfert_dual}
     \begin{aligned}
    &L^*\varphi(x)=\lambda_0 \varphi(x),\, x\in\mathbb{R}_+\\
    &\varphi(x)>0,\, x\in\mathbb{R}_+\\
    &\int_0^\infty \varphi(x)dx=1.
\end{aligned}
\end{equation}
\end{corollary}
\begin{proof}
It suffices to apply the previous approach developed for Proposition \ref{prop:existence_sol_sans_transfert} by replacing $L$ with $L^*$ and considering the operator $Q_\lambda$ defined for $\lambda \in \Sigma$ by
\begin{equation*}
    Q_\lambda\varphi(x):=\frac{\gamma(x)}{D_\lambda(x)}\int_0^\infty \varphi(y)m(x,y)dy.
\end{equation*}
According to Theorem \ref{thm:Jentzsch}, there exists a pair $(\varphi_\lambda, \psi(\lambda)) \in L^1(d\nu) \times \mathbb{R}_+^*$ such that $Q_\lambda\varphi_\lambda = \psi(\lambda)\varphi_\lambda$. 

By the same arguments as in the proof of Proposition \ref{prop:existence_sol_sans_transfert}, there exists a unique parameter $\lambda^* \in \Sigma$ satisfying $\psi(\lambda^*) = 1$. Let $\varphi = \varphi_{\lambda^*}$ denote the associated renormalised eigenfunction. Since $\lambda^* \in \Sigma$, Assumption (H.c) ensures that $\lambda^* > 0$. As previously established, this eigenfunction $\varphi$ is continuous and bounded. 

Because $N \in L^1(dx)$ and $\varphi \in C_b(\mathbb{R}_+)$, their product $N\varphi$ is integrable with respect to the Lebesgue measure. By the duality of the operators $L$ and $L^*$, we can therefore write:
\begin{align*}
    \int_0^\infty N(x)L^*\varphi(x)dx &= \int_0^\infty \lambda^*\varphi(x)N(x)dx, \\
    \int_0^\infty LN(x)\varphi(x)dx &= \lambda^*\int_0^\infty \varphi(x)N(x)dx, \\
    \lambda_0\int_0^\infty \varphi(x)N(x)dx &= \lambda^*\int_0^\infty \varphi(x)N(x)dx.
\end{align*}
Given the strict positivity of $N$ and $\varphi$ almost everywhere, we have $\int_0^\infty \varphi(x)N(x)dx > 0$. We can therefore divide both sides by this integral to conclude that $\lambda^* = \lambda_0$.
\end{proof}

\section{Existence of a solution to the horizontal transfer equation}\label{section:pb_with_transfer}
To prove the existence of a stationary state, our strategy relies on a fixed-point argument, which first requires the careful construction of a dedicated operator. To this end, we embed our target problem into a family of problems parameterized by $\delta \in [0,1]$. This parameter acts as a bridge, connecting a simpler, transfer-free problem ($\delta=0$) to our true, full model ($\delta=1$). For a given function $v$, we will "freeze" the nonlinearity to define a linear eigenvalue problem, whose principal eigenfunction will yield our fixed-point operator.

\textit{Step 1: Truncation and effective transfer.}
Let $\mathcal{A} = \{v \in L^1(d\nu) \mid \varepsilon < \|v\|_\nu < E+1\}$ be the set of admissible profiles, where the constants $\varepsilon, E > 0$ correspond to \textit{a priori} bounds that will be determined later. We denote by $\overline{\mathcal{A}}$ the closure of $\mathcal{A}$.
For any pair $(v, \delta) \in \overline{\mathcal{A}} \times [0,1]$, we first normalize $v$ by setting $\widetilde{v} = |v| / \|v\|_1$. This projection is well-defined because $\|v\|_\nu \geq \varepsilon > 0$ implies that $v$ is not almost everywhere zero, which guarantees $\|v\|_1 > 0$.

Using this normalized profile, we construct an effective non-local transfer term $\widetilde{H}_\delta$. We define $\widetilde{\rho}_\delta = \frac{1}{C} \int_0^\infty r(y) \widetilde{v}(y) \, dy$, and set:
\begin{equation*}
\widetilde{H}_\delta(x) = \delta h(\widetilde{\rho}_\delta) \int_0^\infty \text{sgn}(x-y) \widetilde{v}(y) \, dy.
\end{equation*}
Here, $h(x) = \frac{\tau x}{\beta+\mu x} \mathds{1}_{\{x>0,\,\beta>0\}} + \frac{\tau}{\mu} \mathds{1}_{\{\beta=0\}}$ is a continuous and non-negative function on $\mathbb{R}$. Note that $h$ vanishes on $\mathbb{R}_-$ when $\beta > 0$, but remains a positive constant if $\beta = 0$. The crucial property is that $h$ is always non-negative: this ensures we never invert the sign of $\widetilde{H}_\delta$ (which would unphysically reverse the direction of the transfer) in the event that $\widetilde{\rho}_\delta \leq 0$. Since $\widetilde{v}$ is a probability density and $h$ is bounded, it immediately follows that $\widetilde{H}_\delta$ remains uniformly bounded by a constant $E_\delta$.

\textit{Step 2: The linearized spectral problem.}
With the nonlinearity fixed by $\widetilde{v}$, we consider the following linearized eigenvalue problem for the unknown eigenpair $(\lambda_\delta, \widetilde{u}_\delta)$:
\begin{equation}\label{eq:equation_linearisee_avec_transfert}
\widetilde{u}_\delta(x) \left( \lambda_\delta + \gamma(x) - r(x) - \widetilde{H}_\delta(x) \right) = \int_0^\infty \gamma(y) m(y,x) \widetilde{u}_\delta(y) \, dy, \quad \text{with } \widetilde{u}_\delta > 0.
\end{equation}
To solve this, we introduce the auxiliary operator $T_{\lambda,\delta}$ acting on $L^1(d\nu)$:
\begin{equation}\label{eq:operateur_spectral}
    T_{\lambda,\delta}\phi(x) := \frac{1}{D_{\lambda,\delta}(x)} \int_0^\infty \gamma(y) m(y,x) \phi(y) dy=\frac{1}{D_{\lambda,\delta}(x)} \int_0^\infty\frac{\gamma(y) m(y,x) \phi(y)}{1+|r(y)|}\nu(dy),
\end{equation}

where $D_{\lambda,\delta}(x) = \lambda + \gamma(x) - r(x) - \widetilde{H}_\delta(x)$. This operator is well-defined for any $\lambda$ in the admissible set $\Sigma_\delta = \{\lambda \in \mathbb{R} \mid \forall x \geq 0, \, D_{\lambda,\delta}(x) > 0\}$. 
Because $\widetilde{H}_\delta$ is uniformly bounded, the exact same arguments as in the Proposition \ref{prop:existence_sol_sans_transfert} (Jentzsch theorem, Marek's Min-Max formula, and the Intermediate Value Theorem) guarantee the existence of a unique principal eigenpair $(\lambda_\delta, \widetilde{u}_\delta)$ solving \eqref{eq:equation_linearisee_avec_transfert}.

\textit{Step 3: Definition of the operator $\Gamma$.}
We can now define the mapping $\Gamma: \overline{\mathcal{A}} \times [0,1] \to L^1(d\nu)$ by returning the normalized principal eigenfunction:
\begin{equation}\label{eq:def_fixed_point_operator}
\Gamma(v,\delta) = \frac{\widetilde{u}_\delta}{\|\widetilde{u}_\delta\|_1}.
\end{equation}
This construction directly serves our main objective: for $\delta \in \{0, 1\}$, if $u$ is a fixed point of $\Gamma(\cdot, \delta)$ satisfying $\lambda_\delta > 0$, then $u$ yields a solution to the corresponding stationary problem with $C\rho_\delta = \lambda_\delta$. Specifically, $\delta = 0$ solves equation \eqref{eq:stationnaire_sans_transfert}, while $\delta = 1$ resolves the true, full model \eqref{eq:stationnaire}.
Let us state our main result.
    \begin{theorem}\label{thm:main}
Under assumption \textbf{(H)}, equation \eqref{eq:stationnaire} admits at least one strictly positive solution if the eigenvalue $\lambda_0$ defined in \eqref{eq:pb_sans_transfert} satisfies
\begin{equation}\label{eq:condition_sur_lambda_zero}
\lambda_0 > \frac{\tau}{\mu}\mathds{1}_{\{\beta=0\}}.
\end{equation}
Furthermore, any such solution is continuous and bounded.
\end{theorem}
\begin{remark}
Under assumption (H.c), the principal eigenvalue $\lambda_0$ is guaranteed to be strictly positive. As a result, in the case $\beta > 0$, the existence condition \eqref{eq:condition_sur_lambda_zero} holds without any condition, ensuring the existence of at least one solution to \eqref{eq:stationnaire}.
\end{remark}
\begin{remark}
    Let us demystify \eqref{eq:condition_sur_lambda_zero}. Recall that the transfer intensity is proportional to $I_{TH}(\rho) = \frac{\tau\rho}{\beta + \mu\rho}$ and that the parameter $\rho$ can be viewed as the true population size for the unnormalised equation or the renormalised global fitness of the population. When this parameter is close to 0, it means that the population is near extinction.
    \begin{itemize}
        \item Assuming $\beta,\mu > 0$, as $\rho$ approaches $0$, we observe the linear behaviour $\frac{\tau}{\beta}\rho$. The limiting factor is the encounter time between two bacteria, namely $\beta$. When $\rho$ is sufficiently large, we observe saturation $\frac{\tau}{\mu}$. The limiting factor is the time it takes for the transfer to occur, namely $\mu$. Biologically, one must find a partner (parameter $\beta$) and when there are few individuals, there is less frequent transfer. The transfer itself takes time (parameter $\mu$). Even when the population size is very large, a bacterium cannot receive more than a certain number of plasmids per unit of time.
        \item Assuming $\beta > 0$, and $\mu = 0$, then $I_{TH}(\rho) = \frac{\tau}{\beta}\rho$. The transfer is instantaneous. There is no saturation limit. The only limiting factor is the encounter rate. Note that in our model, since the population size is limited by the competitive term, there is implicitly an inherent saturation.
        \item Assuming $\beta = 0$ and $\mu > 0$, then $I_{TH}(\rho) = \frac{\tau}{\mu}$. The intensity is constant and no longer depends on population density. Finding a partner is instantaneous. The only limiting factor is the transfer time. Consequently, even a small population is subject to transfer. This is the only case where the condition for the existence of a stationary state is $\lambda_0 > \frac{\tau}{\mu}$, and not simply $\lambda_0 > 0$, precisely because the transfer is not proportional to population density, and its deleterious influence is always present.
    \end{itemize}
    \end{remark}
\begin{remark}
All simulations of equation \eqref{eq:integrodiff} where almost-persistent cycles and convergence to the zero state are observed correspond to parameter sets for which equation \eqref{eq:condition_sur_lambda_zero} is not satisfied.
\end{remark}
The proof relies on the application of the Leray-Schauder topological degree theorem, using the parameter $\delta \in [0,1]$ to build a continuous homotopy between the problem without transfer ($\delta = 0$) and the full target problem ($\delta = 1$). Since existence and uniqueness are already established for the uncoupled case $\delta = 0$, the strategy consists in propagating this existence result up to $\delta = 1$. To achieve this, we first establish the continuity and compactness of the fixed-point operator $\Gamma$ on $\overline{\mathcal{A}} \times [0,1]$. We then derive uniform \textit{a priori} estimates for the family of solutions $(u_\delta)_{\delta \in [0,1]}$, which ensures that no fixed point belongs to the boundary $\partial \mathcal{A}$. 
\begin{proof}[Proof of Theorem \ref{thm:main}]

The compactness and continuity of the operator $\Gamma$ defined in \eqref{eq:def_fixed_point_operator} are proved in Proposition \ref{prop:regularite_Gamma}. Let us establish the \textit{a priori} estimates: 

Let $u_\delta$ be a fixed point of $\Gamma(\cdot,\delta)$ in $\mathcal{A}$. We want to establish bounds on $\|u_\delta\|_\nu = \int_0^\infty |u_\delta(x)|(1+|r(x)|)dx$ and we also need to establish bounds on $\lambda_\delta$ to ensure that it is stricly positive, otherwise $\rho_\delta$ would be non-positive and this would correspond to a non-positive population density. Let us rewrite the equation satisfied by $u_\delta$:
\begin{equation}\label{equation_stationnaire_renormalisee_avec_delta}
    \begin{aligned}
        &(r(x)-\gamma(x)-C\rho_\delta)u_\delta(x)+ \int_0^\infty\gamma(y) u_\delta(y)m(y,x)dy+\delta u_\delta(x) \mathcal{H}(u_\delta,\rho_\delta)(x)=0,\\
        &\rho_\delta=\frac{1}{C}\int_0^\infty r(y)u_\delta(y)dy,\\
        &\mathcal{H}(u_\delta,\rho_\delta)(x)= h(\rho_\delta)\int_0^\infty \text{sgn}(x-y)u_\delta(y)dy,\, h(x)=\frac{\tau x}{\beta+\mu x}\mathds{1}_{\{x>0,\,\beta>0\}}+\frac{\tau}{\mu}\mathds{1}_{\{\beta=0\}},\\
        &u_\delta(x)>0\text{ and }\int_0^\infty u_\delta(y)dy=1.
    \end{aligned}
\end{equation}

\underline{Bounds for $\rho_\delta$:}

We begin by establishing uniform bounds for $\rho_\delta$, which subsequently allow us to derive estimates for $\|u\|_\nu$. These norm estimates will, in turn, dictate the choice of the constants $\varepsilon$ and $E$ defining the admissible set $\mathcal{A}$, ultimately enabling us to apply the Leray-Schauder fixed-point theorem.
By the definition of $\rho_\delta$ in \eqref{equation_stationnaire_renormalisee_avec_delta}, we have the immediate upper bound 
    \begin{equation}\label{eq:upper_bound_on_rhodelta}
        \forall\delta\in[0,1],\, \rho_\delta\leq\frac{r_0}{C}.
    \end{equation}
    For the lower bound, let us recall that $\varphi$ is the unique positive solution of \eqref{eq:pb_sans_transfert_dual}. The associated eigenvalue satisfies $\lambda_0 \geq \sup_{x\geq0}(r(x)-\gamma(x))>0$. 
    We have 
    \begin{equation*}
        L^*\varphi=\lambda_0 \varphi=C\rho_0 \varphi
    \end{equation*}
    and 
    \begin{equation*}
        Lu_\delta+\delta u_\delta \mathcal{H}(u_\delta,\rho_\delta)=C\rho_\delta u_\delta.
    \end{equation*}
    Using the duality of $L$ and $L^*$ and the boundedness of $\varphi$, we can write the following identity

    \begin{align*}
    C\rho_\delta\int_0^\infty u_\delta(x)\varphi(x)dx&=\int_0^\infty Lu_\delta(x)\varphi(x)dx+\int_0^\infty \delta u_\delta(x) \mathcal{H}(u_\delta,\rho_\delta)(x) \varphi(x)dx\\
    &=\int_0^\infty u_\delta(x)L^*\varphi(x)dx+\int_0^\infty \delta u_\delta(x) \mathcal{H}(u_\delta,\rho_\delta)(x) \varphi(x)dx \\
    &=C\rho_0\int_0^\infty u_\delta(x)\varphi(x)dx+\int_0^\infty \delta u_\delta(x) \mathcal{H}(u_\delta,\rho_\delta)(x) \varphi(x)dx
    \end{align*}
Given the strict positivity of $\varphi$ and $u_\delta$, we obtain
\begin{equation}\label{eq:debut_priori}
    C(\rho_\delta-\rho_0)=\delta\frac{\int_0^\infty u_\delta(x) \mathcal{H}(u_\delta,\rho_\delta)(x) \varphi(x)dx}{\int_0^\infty u_\delta(x)\varphi(x)dx}.
\end{equation}
We observe that the transfer term is bounded by $\delta h(\rho_\delta)$. Considering the worst case, we get:
\begin{equation*}
    C(\rho_\delta-\rho_0)\geq - \delta h(\rho_\delta)\geq - h(\rho_\delta).
\end{equation*}
We will establish the conditions on $\rho_0,\,\tau,\,\beta$ and $\mu$ for $\rho_\delta$ to be strictly positive.

The non-vanishing of $\rho_\delta$ depends on the parameters $\beta$ and $\mu$. When $\beta > 0$, the estimate reads
\begin{equation*}
C\rho_\delta + \frac{\tau \rho_\delta}{\beta + \mu \rho_\delta} \mathds{1}_{\{\rho_\delta > 0\}} \geq C\rho_0,
\end{equation*}
which immediately implies that $\rho_\delta > 0$ as soon as $\rho_0 > 0$. Conversely, when $\beta = 0$ (and $\mu > 0$), we obtain the bound $C\rho_\delta + \frac{\tau}{\mu} \geq C\rho_0$; in this case, ensuring $\rho_\delta > 0$ requires the threshold condition $\rho_0 > \frac{\tau}{C\mu}$.
\newline

\underline{Bounds for $\|u_\delta\|_\nu$:}

We have the conditions for $\rho_\delta$ to be strictly positive and an upper bound \eqref{eq:upper_bound_on_rhodelta}. We now compute estimates for $\|u_\delta\|_\nu$.

We establish the constants $\varepsilon$ and $E$ defining the open set $\mathcal{A}$, we first consider the lower bound. Since any fixed point $u_\delta$ satisfies \eqref{equation_stationnaire_renormalisee_avec_delta}, we have
\begin{equation*}
\|u_\delta\|_\nu = \int_0^\infty |u_\delta(x)| (1+|r(x)|) \, dx \geq 1 + C\rho_\delta > 1,
\end{equation*}
and we can fix $\varepsilon = \frac{1}{2}$.

For the upper bound, we recall that $u_\delta$ satisfies the fixed-point equation
\begin{equation*}
u_\delta(x) \left( C\rho_\delta + \gamma(x) - r(x) - \delta \mathcal{H}(u_\delta,\rho_\delta)(x) \right) = \int_0^\infty \gamma(y) m(y,x) u_\delta(y) \, dy.
\end{equation*}
By assumption (H.a), there exists $x_0 \geq 0$ such that $r(x_0) = 0$ and $|r(x)| = -r(x)$ for all $x > x_0$. Thus, for $x > x_0$, we can isolate the growth term:
\begin{align*}
|r(x)| u_\delta(x) &= \int_0^\infty \gamma(y) m(y,x) u_\delta(y) \, dy - u_\delta(x) \left( C\rho_\delta + \gamma(x) - \delta \mathcal{H}(u_\delta,\rho_\delta)(x) \right) \\
&\leq \overline{\gamma} \int_0^\infty m(y,x) u_\delta(y) \, dy + u_\delta(x) \left( C\rho_\delta + \overline{\gamma} + \delta \left| \mathcal{H}(u_\delta,\rho_\delta)(x) \right| \right).
\end{align*}
Integrating this inequality over $(x_0, \infty)$ and combining it with the uniform bound \eqref{eq:upper_bound_on_rhodelta}, the boundedness of $\mathcal{H}(u_\delta,\rho_\delta)$, and assumption (H.b) on $m$, we obtain
\begin{align*}
\int_0^\infty |r(x)| u_\delta(x) \, dx &\leq \int_0^{x_0} |r(x)| u_\delta(x) \, dx + \overline{\gamma} \int_0^\infty \int_0^\infty m(y,x) u_\delta(y) \, dy \, dx \\
&\quad + \left( C\rho_\delta + \overline{\gamma} + C(\beta,\mu) \right) \int_{x_0}^\infty u_\delta(x) \, dx \\
&\leq 2r_0 + 2\overline{\gamma} + C(\beta,\mu),
\end{align*}
where the parameter-dependent constant is defined by
\begin{equation*}
C(\beta,\mu) = \frac{\tau}{\mu} \mathds{1}_{\{\beta \geq 0, \, \mu > 0\}} + \frac{\tau r_0}{\beta C} \mathds{1}_{\{\beta > 0, \, \mu = 0\}}.
\end{equation*}
Setting $E = 1 + 2r_0 + 2\overline{\gamma} + C(\beta,\mu)$ leads to  $\|u_\delta\|_\nu \leq E$ for all $\delta$ in $[0,1]$.

In summary, our choice of admissible open set $\mathcal{A}$ is
\begin{equation*}
\mathcal{A} = \left\{ v \in L^1(d\nu) \;\middle|\; \frac{1}{2} < \|v\|_\nu < 2(1+r_0 + \overline{\gamma}) + C(\beta,\mu) \right\}.
\end{equation*}

\underline{Conclusion:}

We are now in a position to apply the Leray-Schauder degree Theorem (see Theorem \ref{thm:leray_schauder} in Appendix \ref{appendix}). To do so, we must verify four fundamental requirements. First, by construction, the admissible set $\mathcal{A} = \{v \in L^1(d\nu) \mid \varepsilon < \|v\|_\nu < E+1\}$ is a bounded open subset of $L^1(d\nu)$. Second, the fixed-point operator $\Gamma: \overline{\mathcal{A}} \times [0,1] \to L^1(d\nu)$ is continuous and compact; the detailed proof of these topological properties is postponed to Section \ref{section:regularity_results}. Third, thanks to the \textit{a priori} estimates established previously, any fixed point $u = \Gamma(u,\delta)$ for $\delta \in [0,1]$ necessarily satisfies $\varepsilon < \|u\|_\nu \leq E$. This ensures that no fixed point lies on the boundary $\partial \mathcal{A}$. Finally, for the transfer-free model ($\delta = 0$), the operator $\Gamma(\cdot, 0)$ admits a unique non-zero fixed point in $\mathcal{A}$, guaranteeing that the Leray-Schauder topological degree is strictly positive.

Having satisfied these conditions, the Leray-Schauder degree Theorem guarantees the existence of at least one fixed point $u = \Gamma(u,1) \in \mathcal{A}$ for the full model. Under condition \eqref{eq:condition_sur_lambda_zero}, this establishes the existence of a solution in $L^1(d\nu)$ to the renormalised stationary equation \eqref{eq:stationnaire}.

Moreover, this solution automatically belongs to $L^1(d\nu)\cap C_b(\mathbb{R}_+)$. Specifically, since $u$ is a fixed point, it explicitly satisfies
\begin{equation}\label{eq:strict_positivity}
        u(x) = \frac{\int_0^\infty\gamma(y) m(y,x) u(y) \, dy}{C\rho + \gamma(x)- r(x) + \mathcal{H}(u,\rho)(x)} > 0.
\end{equation}
Both the numerator and the denominator in this quotient are continuous with respect to $x$. Therefore, we can identify $u$ with its continuous representative.
Its boundedness is a direct consequence of assumption (H.a).
Finally, the strict positivity of $\overline{n}$ is a consequence of \eqref{eq:strict_positivity}.
\end{proof}
\begin{remark}
We once again insist on the condition $\rho_\delta > 0$, which is not the true lower bound of the a priori estimate. Indeed, the fixed-point theorem gives us a solution of mass $1$, but if the corresponding quantity $\rho$ is negative or zero, we are dealing with an extinction. When reconstructing the true solution (which is not necessarily of mass $1$), we obtain a negative solution everywhere, which is neither the sought-after solution nor biologically desirable.
\end{remark}
We conclude this section by investigating the stability of the zero state when $\beta > 0$, which is particularly useful for validating our subsequent numerical results.

\begin{corollary}\label{cor:instabilite_zero_beta_pos}
Assume $\beta>0$. Assumptions \textbf{(H)} ensure that
the principal eigenvalue $\lambda_0$ of the linear transfer-free operator $L$ is strictly positive and it implies that the trivial stationary state $n \equiv 0$ of \eqref{eq:stationnaire} is linearly unstable.
\end{corollary}

\begin{proof}
To investigate the linear stability of the trivial state, we consider a small perturbation $n(t,x) = \varepsilon \phi(t,x)$ with $0 < \varepsilon \ll 1$. The corresponding total mass is $\rho(t) = \varepsilon \int_0^\infty \phi(t,y)dy = \varepsilon \tilde{\rho}(t)$.

Injecting this perturbation into the main integro-differential equation~\eqref{eq:integrodiff}, the global competition term scales as $C \rho n = C \varepsilon^2 \tilde{\rho} \phi = \mathcal{O}(\varepsilon^2)$.

Furthermore, since $\beta > 0$, the denominator of the transfer term remains strictly positive and bounded away from zero as $\rho \to 0$. The horizontal transfer term thus scales quadratically with respect to the perturbation:
\begin{equation*}
    \frac{\tau}{\beta + \mu \varepsilon \tilde{\rho}} (\varepsilon \phi(x)) \int_0^\infty \text{sgn}(x-y) \varepsilon \phi(y) dy = \mathcal{O}(\varepsilon^2).
\end{equation*}

Collecting only the first-order terms in $\varepsilon$, the linearized system around the zero state reduces to:
\begin{equation*}
    \partial_t \phi(t,x) = r(x)\phi(t,x) + \int_0^\infty \left[\gamma(y)\phi(t,y)m(y,x) - \gamma(x)\phi(t,x)m(x,y)\right]dy = L\phi(t,x),
\end{equation*}
where $L$ is the fundamental linear operator without competition and without transfer, given by \eqref{eq:def_linear_operator}. 

The assumption that the transfer-free equation~\eqref{eq:stationnaire_sans_transfert} admits a non-trivial solution with mass $\rho_0 > 0$ implies that the principal eigenvalue $\lambda_0$ of the operator $L$ must be strictly positive (since $\lambda_0 = C \rho_0 > 0$). Consequently, the spectrum of the linearized system around the zero state contains a strictly positive real part, making the trivial state $n \equiv 0$ linearly unstable for any $\beta > 0$.
\end{proof}
\section{Explicit existence conditions and approximations}\label{section:explicit_conditions}
In this section, we adopt simplified assumptions in order to obtain more explicit existence conditions and to highlight the interactions among the model parameters.

\subsection*{Assumption \textbf{(H')}}\label{Assumptions_H_2}
\begin{itemize}
    \item (H'.a) The reproduction rate is given by $r(x) = r_0 - \alpha x^\theta$, with $r_0, \alpha, \theta > 0$. 
    \item (H'.b) The mutation kernel $m$ takes the symmetric, reflected form
    \begin{equation*}
        m(y,x) = \hat{m}(x-y) + \hat{m}(x+y), \quad x,y \geq 0,
    \end{equation*}
    where $\hat{m}$ is a symmetric probability density on $\mathbb{R}$ (typically a Gaussian or Cauchy distribution).
    \item (H'.c) The mutation rate $\gamma > 0$ is constant.
    \item (H.d) $C$ and $\tau$ are positive, $\beta$ and $\mu$ are non negative and non-zero simultaneously.
\end{itemize}
The term $\hat{m}(x-y)$ represents standard unbiased variations from trait $y$ to trait $x$. The second term, $\hat{m}(x+y)$, acts as a reflection boundary condition at the origin $x = 0$. Biologically, any mutant attempting to cross into the non-physical negative trait space ($x < 0$) is folded back into the positive domain $[0, \infty)$, ensuring that the total mutation mass is conserved within the admissible trait space.

\begin{remark}
Compared to assumption \textbf{(H)}, we drop the strong positivity condition \eqref{eq:condition_forte_de_positivite}, meaning that $r_0 > \gamma$ is no longer assumed to hold.
If we maintained the identity \eqref{eq:condition_forte_de_positivite} with the assumption \textbf{(H')}, then we would have
\begin{equation}\label{eq:borne_naive}
    r_0>\lambda_0>r_0-\gamma>0
\end{equation}
which is a simple but too imprecise bound. It is possible to obtain more precise estimates, even when if $\gamma \geq r_0$.
\end{remark}

Under assumption \textbf{(H')}, the operator $L$ is self-adjoint in $L^2$. We use Theorem \ref{thm:min_max_variationnel} which provides the expression
\begin{equation}\label{eq:expression_L2_lambda0}
    \lambda_0:=\sup_{\phi\in L^2(\mathbb{R}_+), \|\phi\|_{L^2}=1}\left(  \int_0^\infty r(x)\phi^2(x)dx-\gamma\left(1-\int_0^\infty\int_0^\infty m(y,x)\phi(x)\phi(y)dydx\right)\right).
\end{equation}
Due to the strict decay of the spectral radius $\lambda \rightarrow \psi(\lambda)$ (proven in Proposition \ref{prop:existence_sol_sans_transfert}), we know that $\lambda_0$ is well-isolated and the largest positive eigenvalue of $L$. Indeed, despite the new hypotheses, the arguments established in Proposition \ref{prop:existence_sol_sans_transfert} remain valid; the only property that no longer holds is the strict positivity of $\lambda_0$.  Because the stationary solution $N$ of \eqref{eq:pb_sans_transfert} is $\nu$-integrable and bounded, it is also an element of $L^2(\mathbb{R}_+)$. Since the operator $L$ is self-adjoint in $L^2$ and bounded above, we can calculate its largest eigenvalue using \eqref{eq:expression_L2_lambda0}.
\begin{proposition}
Under the assumption \textbf{(H')} and if $\hat{m}(x)=\frac{1}{\sqrt{2\pi\sigma^2}} e^{\frac{-x^2}{2\sigma^2}}$, then \eqref{eq:stationnaire} admits a solution if there exists a real $a>0$ such that
\begin{equation}
    r_0>\frac{\alpha}{\sqrt{\pi}}\Gamma\left(\frac{\theta+1}{2}\right)a^{-\tfrac{\theta}{2}}+\gamma\left(1-\sqrt{\frac{2}{2+a\sigma^2}}\right)+\frac{\tau}{\mu}\mathds{1}_{\{\beta=0,\mu>0\}}
\end{equation}
where $\Gamma(x)=\int_0^\infty t^{x-1}e^{-t}dt$. If $\theta=2$ then the condition becomes
\begin{equation}\label{eq:existence_condition_precise_gaussian}
    r_0>\frac{\alpha}{2a}+\gamma\left(1-\sqrt{\frac{2}{2+a\sigma^2}}\right)+\frac{\tau}{\mu}\mathds{1}_{\{\beta=0,\mu>0\}}.
\end{equation}
If moreover, $\sigma$ is close to $0$, we have the approximate condition
\begin{equation}
    r_0>\sigma\sqrt{\frac{\alpha\gamma}{2}}+\frac{\tau}{\mu}\mathds{1}_{\{\beta=0,\mu>0\}}.
\end{equation}
\end{proposition}
\begin{proof}
We calculate the expression given by \eqref{eq:expression_L2_lambda0} with the test function $\phi(x) = e^{-\frac{ax^2}{2}}$, we then find
\begin{align*}
    \lambda_0\geq\lambda(a)&:=\int_0^\infty r(x)e^{-ax^2}dx-\gamma\left( 1-\int_0^\infty\int_0^\infty m(y,x)e^{-a\frac{(x^2+y^2)}{2}}dydx \right)\\
    &=r_0-\frac{\alpha}{\sqrt{\pi}}\Gamma\left(\frac{\theta+1}{2}\right)a^{-\tfrac{\theta}{2}}-\gamma\left(1-\sqrt{\frac{2}{2+a\sigma^2}}\right).
\end{align*}
By studying the variations of $\lambda$, one can determine the maximum but the expression is not explicit. Therefore, one proceeds via a Taylor expansion for $\sigma$ in 0 (in the case $\theta=2$) :
\begin{equation*}
    \lambda(a)=\frac{\alpha}{2a}+\frac{\gamma\sigma^2}{4}a+o(\sigma^2).
\end{equation*}
It is then sufficient to find the maximum for this expression.
\end{proof}
\begin{remark}
If the average encounter time between two individuals is non-zero (i.e., $\beta > 0$), we deduce that the population survives if the maximum reproduction rate exceeds the combined influence of selection ($\alpha$), mutation frequency ($\gamma$), and mutation variance ($\sigma^2$). If the parameter $\beta$ is zero, it must be that the maximum reproduction rate is greater than the combined influence of selection, mutations, and transfer intensity $(\tau/\mu)$.
\end{remark}
\begin{remark}
Assuming that $r_0 = \gamma$, then \eqref{eq:existence_condition_precise_gaussian} can be written as 
\begin{equation*}
    \sqrt{\frac{2}{2+a\sigma^2}}>\frac{\alpha}{2a}.
\end{equation*}
We can easily find values of $\sigma$ and $\alpha$ such that there exists an $a>0$ which satisfies this inequality. In this regard, \eqref{eq:existence_condition_precise_gaussian} allows one to circumvent \eqref{eq:condition_forte_de_positivite} and \eqref{eq:borne_naive}.
\end{remark}
Using a similar reasoning, we obtain a condition in the case of Cauchy mutations.
\begin{corollary}
Under assumption \textbf{(H')}, $\theta=2$ and $\hat{m}(x)=\frac{\sigma}{\pi(x^2+\sigma^2)}$, \eqref{eq:stationnaire} admits a solution if there exists a real $a>0$ such that
\begin{equation*}
    r_0>\frac{\alpha}{2a}+\gamma\left(1-e^{a\frac{\sigma^2}{4}}\text{erfc}\left(\frac{\sigma\sqrt{a}}{2}\right)\right)+\frac{\tau}{\mu}\mathds{1}_{\{\beta=0,\mu>0\}}
\end{equation*}
where
\begin{equation*}
    \text{erfc}(x)=\frac{2}{\sqrt{\pi}}\int_x^\infty e^{-t^2}dt.
\end{equation*}
When $\sigma$ is close to $0$, we have the approximate condition
    \begin{equation}
        r_0>\frac{\alpha}{2}F^2+\frac{\gamma\sigma}{\sqrt{\pi}}F+\frac{\tau}{\mu}\mathds{1}_{\{\beta=0,\mu>0\}},\,\text{where }F=\left(\frac{\alpha\sqrt{\pi}}{\gamma\sigma}\right)^{1/3}.
    \end{equation}
\end{corollary}
\section{Proofs of the regularity of the operators}\label{section:regularity_results}
In this section, we prove the compactness of the operators $T_\lambda$ and $\Gamma$.

\begin{lemma}\label{lemma:compactness_of_operator_without_transfer}
Recall that $D_\lambda(x)=\lambda+\gamma(x)-r(x)$ and $\Sigma=\{\lambda\in\mathbb{R}|\,\forall x\geq0,\, D_\lambda(x)>0\}$. Under assumption \textbf{(H)} or \textbf{(H')}, for all $\lambda\in\Sigma$, the operator $T_\lambda : L^1(d\nu)\rightarrow L^1(d\nu)$ defined by
\begin{equation*}
    T_\lambda\phi(x):=\frac{1}{D_\lambda(x)}\int_0^\infty \gamma(y)m(y,x)\phi(y)dy,\, \forall x\geq0
\end{equation*}
is compact.
\end{lemma}
\begin{proof}
Let $B_1 = \{f \in L^1(d\nu) \mid \|f\|_\nu \leq 1\}$ be the unit ball. We apply Theorem \ref{thm:kolmogorov} in $L^1(d\nu)$ (which is possible since the $L^1(d\nu)$ and $L^1(dx)$ are isometric).

Let us begin by showing that the operator is bounded, i.e.,
\begin{equation*}
    \exists C > 0, \forall f \in B_1, \int_0^\infty |T_\lambda f(x)|(1+|r(x)|)dx \leq C.
\end{equation*}
Let $f$ be an element of $B_1$. We have
\begin{equation*}
        \int_0^\infty |T_\lambda f(x)|(1+|r(x)|)dx \leq \overline{\gamma} \int_0^\infty |f(y)| \left(\int_0^\infty \frac{1+|r(x)|}{D_\lambda(x)}m(y,x)dx\right)dy.
\end{equation*}
The term $\frac{1+|r(x)|}{D_\lambda(x)}$ is bounded uniformly in $x$, and $\int_0^\infty m(y,x)dx = 1$, so there exists a $C > 0$ such that
\begin{equation*}
    \int_0^\infty |T_\lambda f(x)|(1+|r(x)|)dx \leq C \int_0^\infty |f(y)|dy \leq C.
\end{equation*}

Now, let us demonstrate the decay at infinity, i.e.,
\begin{equation*}
    \forall \varepsilon > 0, \exists R > 0, \forall f \in B_1, \int_R^\infty |T_\lambda f(x)|(1+|r(x)|)dx < \varepsilon.
\end{equation*}
We define
\begin{equation*}
    I_R = \overline{\gamma} \int_R^\infty \frac{1+|r(x)|}{D_\lambda(x)} \left(\int_0^\infty m(y,x)|f(y)|dy \right)dx.
\end{equation*}
We know that there exists a $C > 0$ which bounds the fraction $\frac{1+|r(x)|}{D_\lambda(x)}$. By Fubini's theorem, we have

    \begin{align*}
        I_R &\leq \overline{\gamma} C \int_R^\infty \int_0^\infty m(y,x)|f(y)|dy dx = \overline{\gamma} C \int_0^\infty |f(y)| \left(\int_R^\infty m(y,x)dx\right) dy.
    \end{align*}

We denote
\begin{equation*}
    S_R(y) = \int_R^\infty m(y,x)dx.
\end{equation*}
We know that, for a fixed $y$, $S_R(y)$ tends to $0$ as $R$ tends to infinity. We split the integral:
\begin{equation*}
    I_R \leq \overline{\gamma} C \left(\int_0^A |f(y)| S_R(y) dy + \int_A^\infty |f(y)| S_R(y) dy\right).
\end{equation*}
We know that $\|S_R\|_\infty \leq 1$, $\|f\|_\nu \leq 1$, and the function $x \mapsto |r(x)|$ is increasing after a certain threshold. For a sufficiently large $A$, we obtain

\begin{align*}
    \int_A^\infty |f(y)| S_R(y) dy &\leq \frac{1}{1+|r(A)|} \int_A |f(y)| (1+|r(y)|) dy \\
    &\leq \frac{1}{1+|r(A)|} \|f\|_\nu \\
    &< \frac{\varepsilon}{2} \text{ (for } A \text{ sufficiently large).}
\end{align*}

Now that $A$ is fixed, let us consider the other integral term. Let $y$ be in $[0,A]$:
\begin{equation*}
    \int_0^A |f(y)| S_R(y) dy \leq \left( \sup_{y\in[0,A]} S_R(y) \right) \|f\|_\nu < \frac{\varepsilon}{2} \text{ for } R \text{ sufficiently large.}
\end{equation*}
By virtue of the continuity of $y\mapsto S_R(y)$, the supremum is attained. For $R$ sufficiently large, we have $I_R < \varepsilon$.
It remains to show property \eqref{eq:continuite_translations}. Thanks to the uniform decay at infinity established above, it is sufficient to prove that the family of functions $\{T_\lambda f \mid f \in B_1\}$ is uniformly equicontinuous on any compact interval $[0,R]$.
We apply Theorem \ref{thm:arzela_ascoli}. We want to show that 
\begin{equation*}
    \forall \varepsilon > 0, \exists \delta > 0, \forall x, y \in [0, R], \forall f \in B_1, |x-y| < \delta \implies |T_\lambda f(x) - T_\lambda f(y)| < \varepsilon.
\end{equation*}
We then use the bound $\overline{M}$ defined in (H.b):

    \begin{align*}
        |T_\lambda f(x)-T_\lambda f(y)|&\leq \frac{\overline{\gamma}}{D_\lambda(x)}\left( \sup_{z\geq0} |m(z,x)-m(z,y)|\right)\\
        &+2\frac{\overline{\gamma}\overline{M}}{D_\lambda(x)D_\lambda(y)}(|r(x)-r(y)|+|\gamma(x)-\gamma(y)|).
    \end{align*}

By the uniform continuity of $m$ and the continuity of $r$ and $\gamma$, we deduce the desired result.
We conclude that the auxiliary operator $T_\lambda$ is compact.
\end{proof}
\begin{proposition}\label{prop:regularite_Gamma}
The operator $\Gamma$ defined in Section \ref{section:pb_with_transfer} is compact.
\end{proposition}
\begin{proof}
We are working on $L^1(d\nu)$. Let $\mathcal{F} = \overline{\mathcal{A}} \times [0,1]$ be the domain of definition of $\Gamma$ which is closed and bounded.

We start by proving the continuity of $\Gamma$.

Let $(v_n, \delta_n)$ be in $\mathcal{F}$ which converges to $(v, \delta)$ (in $\mathcal{F}$). We must show that the sequence of images $(\Gamma(v_n, \delta_n))$ converges to $\Gamma(v,\delta)$ in $L^1(d\nu)$. We recall the notations:

    \begin{align*}
        \widetilde{v} := \frac{|v|}{\|v\|_1}, &\qquad \widetilde{v}_n := \frac{|v_n|}{\|v_n\|_1},\\
        \widetilde{\rho} := \frac{1}{C} \int_0^\infty \widetilde{v}(x)(1+|r(x)|)dx, &\qquad\widetilde{\rho}_n := \frac{1}{C} \int_0^\infty \widetilde{v}_n(x)(1+|r(x)|)dx,\\
        H(x) := \delta h(\widetilde{\rho}) \int_0^\infty \text{sgn}(x-y)\widetilde{v}(y)dy,&\qquad H_n(x) := \delta_n h(\widetilde{\rho}_n) \int_0^\infty \text{sgn}(x-y)\widetilde{v}_n(y)dy, \\ \\
    \end{align*}
Since $v \in \overline{\mathcal{A}}$, we have $\|v\|_\nu \geq \varepsilon > 0$, which implies that the mapping $v \mapsto \widetilde{v} := \frac{|v|}{\|v\|_1}$ is continuous. Thus, $\widetilde{v}_n$ converges to $\widetilde{v}$ in $L^1(d\nu)$, a result that can be justified via dominated convergence.

As a consequence, we estimate
\[
|\widetilde{\rho}_n - \widetilde{\rho}| = \frac{1}{C} \left| \int_0^\infty(1+ |r(x)|)(\widetilde{v}_n(x) - \widetilde{v}(x))\,dx \right| \leq \frac{1}{C} \|\widetilde{v}_n - \widetilde{v}\|_\nu,
\]
which proves the convergence of $\widetilde{\rho}_n$ to $\widetilde{\rho}$. By continuity of $h$, the coefficient $h(\widetilde{\rho}_n)$ also converges to $h(\widetilde{\rho})$. 

Furthermore, because $\|\widetilde{v}_n - \widetilde{v}\|_1 \to 0$, the integral $\int_0^\infty \text{sgn}(x-y)\widetilde{v}_n(y)\,dy$ converges uniformly in $x$ to $\int_0^\infty \text{sgn}(x-y)\widetilde{v}(y)\,dy$. Combining these convergence results, we conclude that $H_n(x) := \delta_n h(\widetilde{\rho}_n) \int_0^\infty \text{sgn}(x-y)\widetilde{v}_n(y)\,dy$ converges uniformly on $\mathbb{R}_+$ to $H(x) := \delta h(\widetilde{\rho}) \int_0^\infty \text{sgn}(x-y)\widetilde{v}(y)\,dy$.

We now turn to the family of operators $(T^{(n)}_\lambda)_{n,\lambda}$ defined below. For all $n$, let us set
\begin{align*}
    S_n &:= \sup_{x \geq 0} \{r(x) - \gamma(x) + H_n(x)\}, \\
    S   &:= \sup_{x \geq 0} \{r(x) - \gamma(x) + H(x)\}.
\end{align*}
Recalling the notations defined in Section \ref{section:pb_with_transfer}, we emphasize that the quantity $\inf \Sigma_\delta$ defined in \eqref{eq:operateur_spectral} is equal to $(S,+\infty)$. 
Given that we are now dealing with a sequence of admissible sets indexed by $n$, proving the convergence of this sequence is a necessary preliminary step to establish the convergence of all other quantities indexed by $n$.
Since $H_n$ converges uniformly on $\mathbb{R}_+$ and $r$ is decreasing, the supremum is attained on a compact set, leading $S_n$ to converges towards $S$.

For each $\lambda \in (S_n, \infty)$, we consider the operator
\[
T_\lambda^{(n)}\phi(x) := \frac{1}{D^{(n)}_\lambda(x)} \int_0^\infty \gamma(y) m(y,x) \phi(y)\,dy,
\]
where $$D^{(n)}_\lambda (x) := \lambda + \gamma(x) - r(x) - H_n(x).$$

Let $\psi^{(n)}(\lambda)$ be the spectral radius of the operator. Following the same arguments as in the case without transfer, there exists a unique $\lambda^*_n$ such that $\psi^{(n)}(\lambda^*_n) = 1$ for each $n$, and similarly a unique solution $\lambda^*$ to $\psi(\lambda^*) = 1$ for the corresponding limit $(v,\delta)$.

Because the spectral radius diverges as $\lambda$ approaches $S$, we have $\lambda^* > S$. Hence, there exists $\eta > 0$ such that $\lambda^* \geq S + 2\eta$. For $n$ sufficiently large, we have $S_n < S + \eta$, which allows us to define the compact interval $K = [\lambda^* - \eta, \lambda^* + \eta] \subset \bigcap_{n \geq N} (S_n, +\infty)$. 

On this interval $K$, the term $D^{(n)}_\lambda$ converges to $D_{\lambda,\delta}$ and remains bounded below by a strictly positive constant. This ensures the convergence in operator norm:
\begin{equation*}
    \sup_{\lambda \in K} \sup_{\|f\|_\nu \leq 1} \|T_{\lambda}^{(n)}f - T_\lambda f\|_\nu \longrightarrow 0 \quad \text{as } n \to \infty.
\end{equation*}

By continuity of the spectral radius with respect to the operator norm, $\psi_n(\lambda) \to \psi(\lambda)$ pointwise on $K$. Since the sequence $(\psi_n)$ consists of decreasing functions converging to a continuous limit, Dini's Theorem \eqref{thm:dini} yields uniform convergence of $(\psi_n)_n$ on $K$. Consequently, the condition $\psi^{(n)}(\lambda^*_n) = 1 = \psi(\lambda^*)$ implies that $\lambda^*_n \to \lambda^*$.

Finally, applying Theorem \ref{thm:continuite_du_rayon_spectral}, the normalised eigenvector $u_n$ associated with $\lambda^*_n$ converges in $L^1(d\nu)$ to the eigenvector $u$ associated with $\lambda^*$. 

Writing $\Gamma(v_n,\delta_n) = u_n := \frac{u^*_n}{\|u^*_n\|_1}$, the convergence $u^*_n \to u^*$ together with $\|u^*\|_1 > 0$ ensures that $u_n \to u$ in $L^1(d\nu)$, thereby establishing the overall continuity of the operator $\Gamma$.
\newline

We now prove the compactness of the operator $\Gamma$. Given a bounded sequence $(v_n, \delta_n)$ in $\mathcal{F}$, let $u_n = \Gamma(v_n, \delta_n)$. To show that the sequence $(u_n)$ is relatively compact, we rely on Theorem \eqref{thm:kolmogorov}. The approach is structured around three key points: uniform boundedness, uniform tightness, and equicontinuity.

To establish the uniform boundedness, we first note that, by construction, $\|u_n\|_1 = 1$. The equation reads:
\begin{equation*}
    r(x)u_n(x)=(\lambda^*_n+\gamma(x)-H_n(x))u_n(x)-\int_0^\infty\gamma(y) m(y,x)u_n(y)dy.
\end{equation*}
Introducing the notation $D_n(x) := \lambda^*_n+\gamma(x)-H_n(x)-r(x)$, and knowing that
\begin{equation*}
    |H_n(x)|\leq \overline{H}:=\frac{\tau}{\mu}\mathds{1}_{\{\beta=0,\mu>0\}}+(E+1)\frac{\tau}{\beta}\mathds{1}_{\{\beta>0\}},
\end{equation*}
integrating the initial equation over the whole domain gives $\lambda^*_n\leq 2(E+1)+\overline{\gamma}$.
Since $\lambda^*_n>S_n$, there exists a lower bound $L_{min}>-\infty$ such that $\lambda^*_n\in [L_{min},L_{max}]$ for all $n$. We then set $\Lambda=\max(|L_{min}|,L_{max})$. 

Furthermore, since the constant $A$ is chosen such that $r(A) = -r_0$, we have $\int_0^A |r(x)|u_n(x)dx \leq r_0$. For $x > A$, the fact that $|r(x)| = -r(x)$ allows us to write:
\begin{equation*}
    D_n(x) = \lambda^*_n + |r(x)| + \gamma(x) - H_n(x).
\end{equation*}
Since $D_n > 0$, the triangle inequality yields:
\begin{equation*}
    |r(x)| \leq D_n(x) + \Lambda + \overline{\gamma} + \overline{H}.
\end{equation*}
Multiplying this expression by $u_n$ and integrating, while using the identity $D_n(x)u_n(x) = \int_0^\infty\gamma(y) m(y,x)u_n(y)dy$, we get:

    \begin{align*}
        u_n(x)|r(x)| &\leq D_n(x)u_n(x) + (\Lambda + \overline{\gamma} + \overline{H})u_n(x) \\
        &\leq \overline{\gamma}\int_0^\infty m(y,x)u_n(y)dy + (\Lambda + \overline{\gamma} + \overline{H})u_n(x).
    \end{align*}

Integration over the global domain finally provides:
\begin{equation*}
    \int_0^\infty |r(x)|u_n(x)dx \leq r_0 + 2\overline{\gamma} + \Lambda + \overline{H}.
\end{equation*}
Because this bound is independent of $n$ and $x$, the sequence $(u_n)$ is uniformly bounded in the space $L^1(d\nu)$.

We must then prove the uniform tightness of the sequence, which means verifying that the mass does not escape to infinity. For a sufficiently large $A$ and $x > A$, we have the inequality:
\begin{equation*}
    u_n(x)|r(x)| \leq \overline{\gamma} \int_0^\infty m(y,x)u_n(y)dy + Lu_n(x),
\end{equation*}
where $L$ is a constant independent of $n$. Fubini's theorem allows us to deduce:
\begin{equation*}
    \int_A^\infty u_n(x)|r(x)|dx \leq \overline{\gamma} \int_0^\infty u_n(y)\left(\int_A^\infty m(y,x)dx\right)dy + L \int_A^\infty u_n(y)dy.
\end{equation*}
Let $S_A(y) = \int_A^\infty m(y,x)dx$. The second term on the right-hand side can be made arbitrarily small for a sufficiently large choice of $A$ because by construction, $\|u_n\|_1=1$. We now focus on the first term. The function $S_A(y)$ tends to $0$ as $A$ tends to infinity for a fixed $y$, and $S_A \leq 1$. Thanks to the uniform bound established previously, there exists a constant $L' > 0$ such that $\sup_n \int_0^\infty u_n(x)|r(x)|dx \leq L'$. 

By splitting the integral at a point $x_0$, we first evaluate the part on $[x_0, \infty)$:
\begin{equation*}
    \int_{x_0}^\infty u_n(y)S_A(y)dy \leq \int_{x_0}^\infty u_n(y)dy \leq \frac{1}{\min_{z \geq x_0}|r(z)|} \int_{x_0}^\infty u_n(y)|r(y)|dy \leq \frac{L'}{\min_{z \geq x_0}|r(z)|}.
\end{equation*}
Since the function $|r|$ grows without bound towards infinity, we can fix a sufficiently large $x_0$, independent of $n$ and $A$, so that the bound is strictly less than a given $\varepsilon > 0$. On the compact segment $[0, x_0]$, the continuity of $m$ ensures that the mass sent to infinity, bounded by $\sup_{y \in [0, x_0]} S_A(y)$, tends uniformly to $0$ as $A$ increases. We can therefore choose a sufficiently large $A$ to make this contribution also less than $\varepsilon$.
By further increasing $A$ if necessary so that $L \int_A^\infty u_n(y) dy < \varepsilon$, we obtain
\begin{equation*}
    \int_A^\infty u_n(x)|r(x)|dx \leq 3\varepsilon.
\end{equation*}
The sequence $(u_n)$ is consequently uniformly tight in $L^1(d\nu)$.

The last step is to verify the $L^1(d\nu)$-equicontinuity. Thanks to the uniform tightness established above, it suffices to prove that the sequence $(u_n)$ is uniformly equicontinuous on any compact subset of $\mathbb{R}_+$. The functions satisfy:
\begin{equation*}
    u_n(x) = \frac{1}{D_n(x)} \int_0^\infty \gamma(y)m(y,x) u_n(y) dy.
\end{equation*}
Given the uniform $L^1$-bound on $(u_n)$, the uniform continuity of $m$, and the local equicontinuity of $D_n$, this local uniform equicontinuity immediately follows. Combined with the uniform tightness, this fulfills the Kolmogorov-Riesz-Fréchet criteria. The sequence $(u_n)$ is therefore relatively compact.
\end{proof}
\section{Numerical Study}\label{sec:numerical_study}

To complement the analytical existence results, we numerically investigate the convergence dynamics toward equilibrium and assess the linear stability of the stationary branches. Indeed, while the existence of a stationary profile has been rigorously established in the previous sections, proving the uniqueness of the non-trivial solution and its global asymptotic stability remains an open challenge. To the best of our knowledge, standard analytical tools such as General Relative Entropy methods (see \cite{michel2007general, michel2005general, perthame2007transport}) cannot be straightforwardly applied here. The fundamentally different natures of the linear variation operator and the non-local, non-linear transfer term preclude the derivation of a common Lyapunov functional, and other analytical strategies we explored proved equally unsuccessful in overcoming this structural hurdle. A promising theoretical avenue lies in shifting to the cumulative formulation of the integro-differential equation. In this framework, not only does the transfer term remarkably acquire a definite negative sign, potentially unlocking new entropy-based approaches, but its inherent structure also naturally suggests the use of Optimal Transport theory (exploiting 1D transport maps via inverse cumulative distributions) to establish uniqueness. Nevertheless, to bypass this theoretical bottleneck in the present work, we resort to a robust numerical framework. The three parameters of interest are $\tau,\beta$ and $\mu$.

First, the time-dependent integro-differential problem is integrated using a forward Euler scheme, combining left-sided finite differences with left Riemann sums for the non-local integral terms. This explicit time-stepping captures the long-time asymptotic behaviour of the population—allowing us to numerically address the open question of convergence—and provides a robust initial guess of the stationary profiles. In a second phase, this preliminary profile is refined to machine precision using a Newton-Raphson solver applied directly to the discretized stationary system. 

Crucially, standard numerical evaluation of the stationary system encounters a major obstacle at the frequency-dependent limit ($\beta = 0$). When assessing the stability of the trivial zero state, the social interaction term involves a division by the total mass $\rho$, leading to a singularity that immediately causes the Jacobian matrix to blow up. It is precisely here that the renormalised formulation introduced in \eqref{eq:stationnaire} proves its numerical power. By evaluating the system through the decoupled pair $(u, \rho)$, where the density is expressed as $\overline{n}(x) = \rho u(x)$ with $\int u(x)dx = 1$, the mass parameter $\rho$ simplifies within the transfer integral. This regularization eliminates the singularity entirely, providing a well-defined Jacobian that allows our solver to reliably track solution branches and compute the spectrum even as the population mass vanishes. 

Finally, while the time-stepping scheme converges robustly for $\beta > 0$, the relaxation dynamics suffer from severe critical slowing down as $\beta \to 0$. To bypass these prohibitive computation times, we access the singular limit $\beta=0$ via numerical parameter continuation from strictly positive $\beta$ values. Furthermore, to investigate the parameter $\tau$ and navigate potential turning points or bifurcation thresholds where standard Newton iterations diverge, we augment our solver with a pseudo-arclength continuation scheme. This parametrisation ensures robust tracking of the stationary branches across folds and unstable regimes.
For a reference on the algorithms used, particularly for the stability analysis, see \cite{kuznetsov1998elements}.

\subsection{Convergence to the equilibrium}
\begin{figure}[h!]
    \centering
    \begin{subfigure}[t]{0.55\textwidth}
        \centering
        \includegraphics[width=\linewidth, height=4.2cm, keepaspectratio]{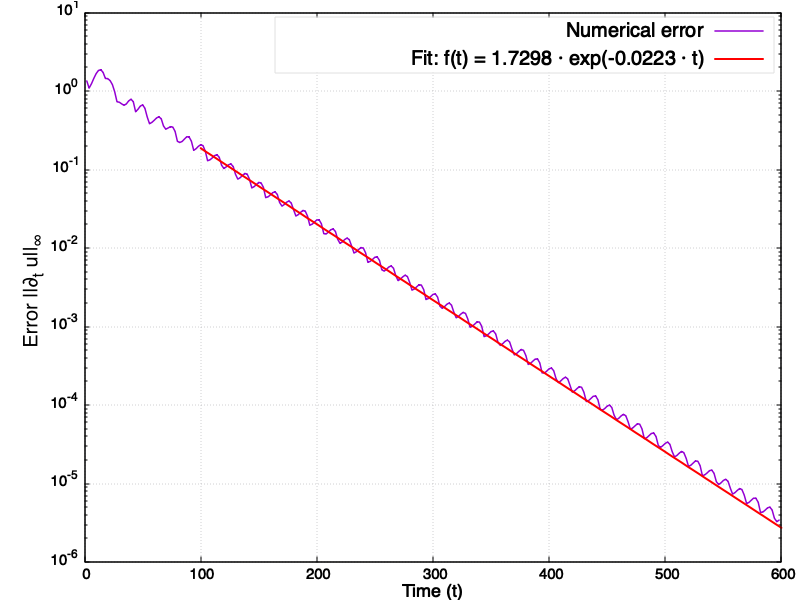}
        \caption{Time evolution of the error $\|\partial_t n(t,\cdot)\|_\infty$ for $\tau = 1,\, \beta = 2$.}
        \label{fig:error_decay}
    \end{subfigure}
    \hfill
    \begin{subfigure}[t]{0.4\textwidth}
        \centering
        \includegraphics[width=1.15\linewidth, height=5.5cm, keepaspectratio]{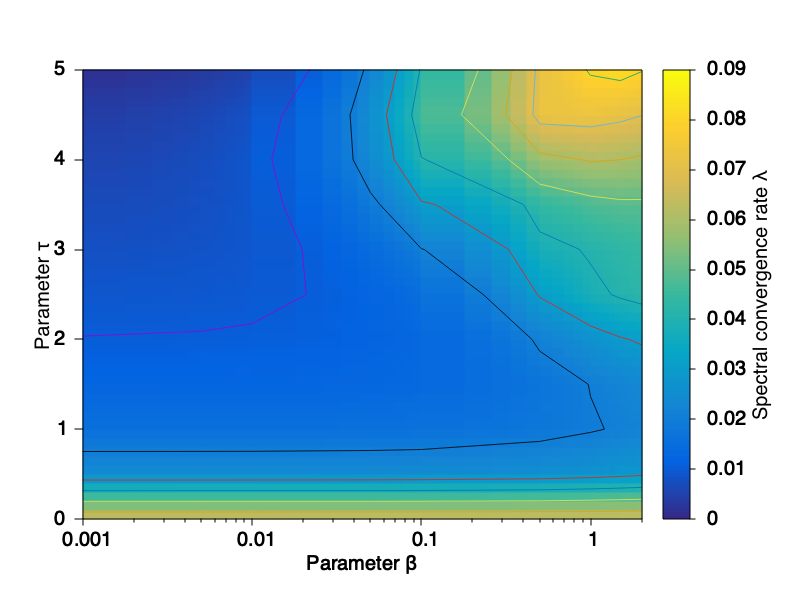}
        \caption{Convergence rate $\lambda$ as a function of $(\beta, \tau)$ in logarithmic scale - $\mu=1$.}
        \label{fig:rate_surface}
    \end{subfigure}
    \caption{Exponential convergence toward equilibrium and parameter sensitivity of the decay rate $\lambda$.}
    \label{fig:convergence_study}
\end{figure}

Across all simulations, the solution systematically converges toward an equilibrium, either the positive stationary state or the zero state. This convergence is noticeably faster when $\beta > 0$.
In order to characterise this convergence to the equilibrium, we simulated the integro-differential equation starting from a uniform initial condition for various values of $\tau$ and $\beta$. Plotting the error on a logarithmic scale confirms that the convergence becomes exponential after an initial transient regime. Specifically,
\begin{equation*}
    \|\partial_t n(t,\cdot)\|_\infty \sim A e^{-\lambda t},
\end{equation*}
where the decay rate $\lambda$ corresponds to the real part of the leading eigenvalue of the linearized system around the equilibrium. Larger values of $\lambda$ reflect faster relaxation toward the stationary state following a perturbation. Here, $\lambda$ is estimated via a linear least-squares fit on the log-transformed error over time, which provides a faster and simpler alternative to extracting the Jacobian spectrum. Figure~\ref{fig:convergence_study} illustrates this exponential decay alongside the dependence of $\lambda$ on $(\beta, \tau)$, where a logarithmic scale is used to capture the singular behaviour as $\beta \to 0$.

As illustrated in Figure~\ref{fig:error_decay}, $\|\partial_t n(t,\cdot)\|_\infty$ exhibits a linear decrease on a logarithmic scale, validating the exponential convergence approximation $A e^{-\lambda t}$. 
Figure~\ref{fig:rate_surface} displays the dependency of this convergence rate $\lambda$ on the parameters $(\beta, \tau)$, highlighting a severe drop in $\lambda$ as $\beta$ approaches $0$. 
This critical slowing down explains why standard time-stepping schemes struggle near the singular limit.

\begin{figure}[h!]
    \centering
    \includegraphics[width=0.48\textwidth]{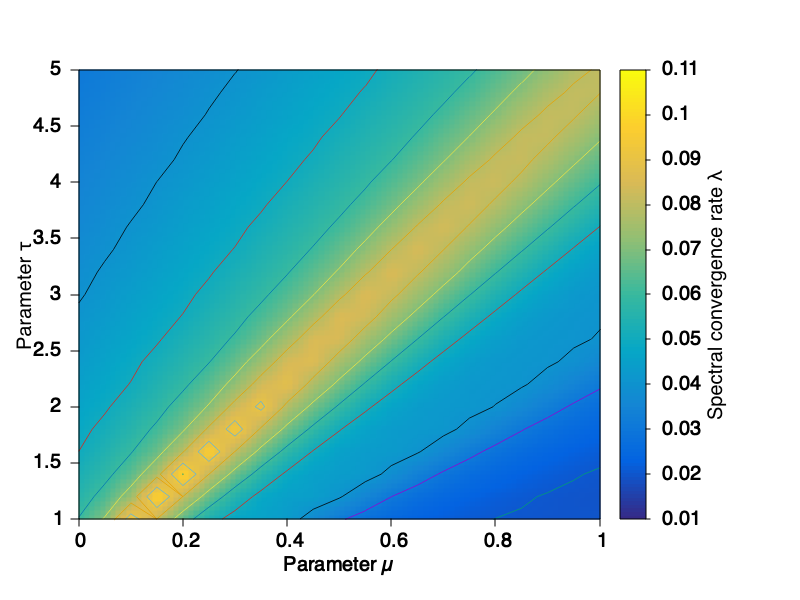}
    \caption{Heatmap of the spectral convergence rate $\lambda$ in the $(\mu, \tau)$ plane (for $\beta = 1$), highlighting a diagonal ridge of optimal relaxation speed.}
    \label{fig:convergence_mu_tau}
\end{figure}

Furthermore, investigating the convergence rate within the $(\mu, \tau)$ parameter space (at a fixed $\beta = 1$) reveals a highly non-monotonic relaxation behaviour. As depicted in Figure~\ref{fig:convergence_mu_tau}, the spectral decay rate $\lambda$ is maximized along a distinct diagonal ridge. This indicates an optimal balance between the social transfer intensity $\tau$ and the mass-dependent saturation parameter $\mu$. Deviating from this optimal ratio—either by independently increasing the transfer strength or the saturation penalty—severely hinders the system's ability to dissipate perturbations, driving the convergence rate down by a full order of magnitude.

\subsection{Study of the stability}
Analysing the stability of the positive steady-state solution is not as straightforward as that of the zero solution in Corollary \ref{cor:instabilite_zero_beta_pos}. This is the purpose of this section. As with the existence criterion, a change in behaviour is observed whether $\beta$ is strictly positive or zero.
\subsubsection{When $\beta>0$}
\begin{figure}[!htbp]
    \centering
    \begin{subfigure}[t]{0.48\textwidth}
        \centering
        \includegraphics[width=\linewidth]{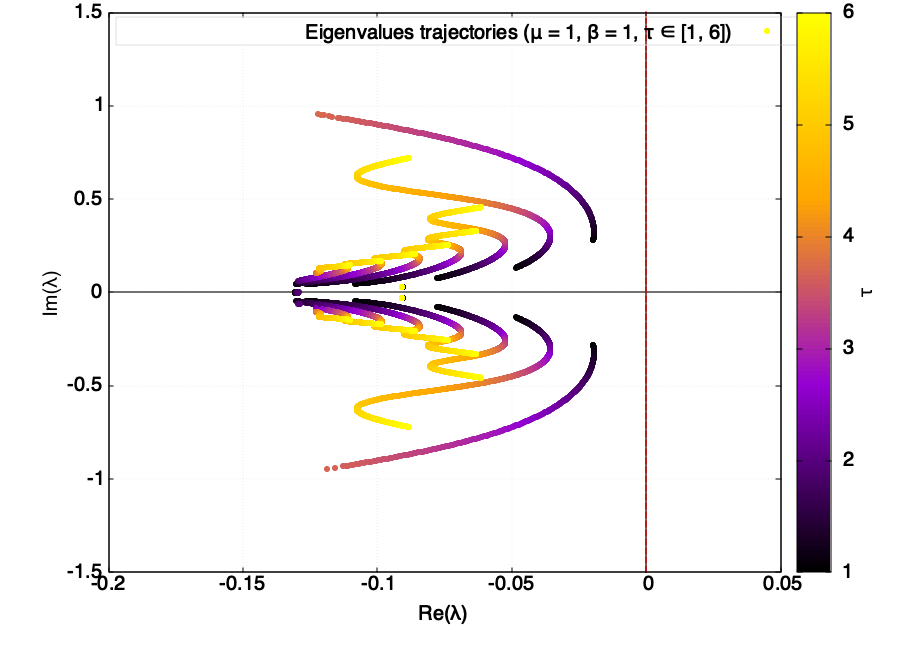}
        \caption{When $\beta=\mu=1$.}
        \label{fig:dominant_eigenvalue_beta>0}
    \end{subfigure}\hfill
    \begin{subfigure}[t]{0.48\textwidth}
        \centering
        \includegraphics[width=\linewidth]{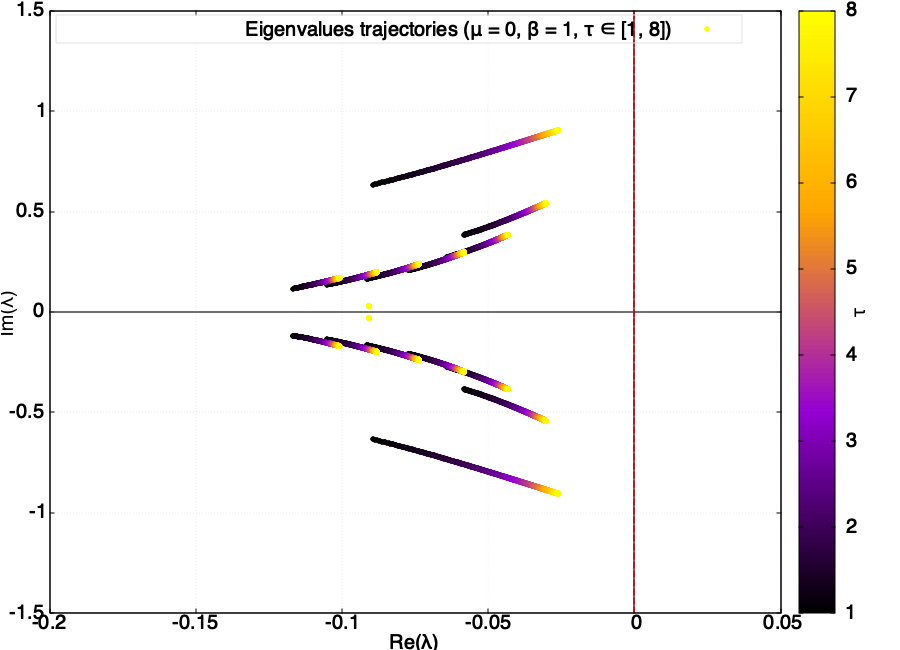}
        \caption{When $\beta=1,\,\mu=0$.}
        \label{fig:jacobian_spectrum_beta>0}
    \end{subfigure}
    \vspace{-2mm}
    \caption{Spectrum of the Jacobian in $\mathbb{C}$ as a function of $\tau$. For each $\tau$ corresponding to a colour, the spectrum is a discrete collection of conjugated points. How the eigenvalues evolve with $\tau$ is seen by following the gradient of colours.}
    \label{fig:stability_analysis_beta>0}
    
    \vspace{3mm}
    \includegraphics[width=0.45\textwidth]{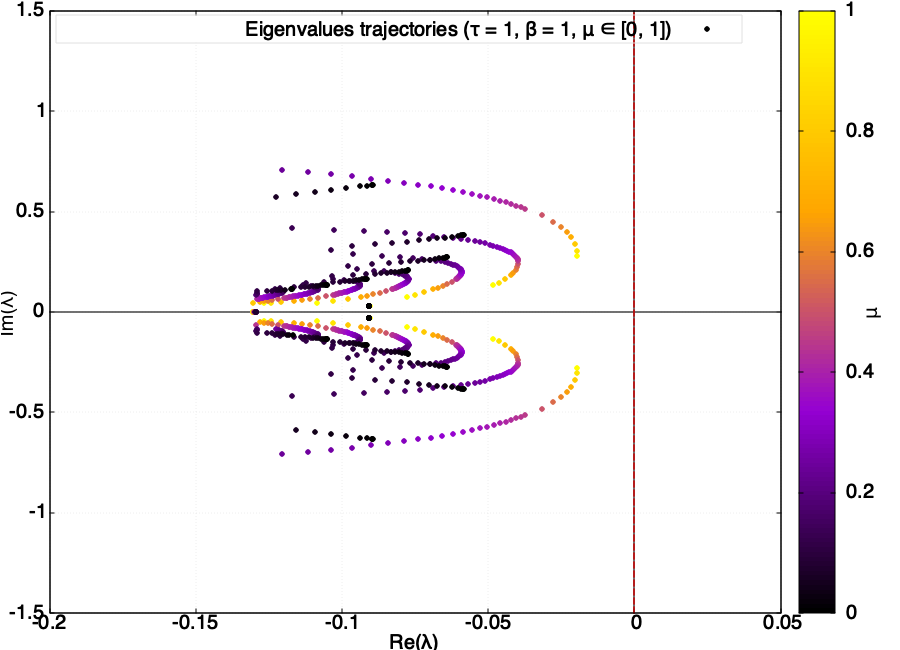}
    \vspace{-2mm}
    \captionof{figure}{Spectrum of the Jacobian in $\mathbb{C}$ as a function of $\mu$ when $\tau=1,\,\beta=1$.}
    \label{fig:stability_analysis_mu_to_0}

    \vspace{3mm}
    \begin{subfigure}[t]{0.48\textwidth}
        \centering
        \includegraphics[width=\linewidth]{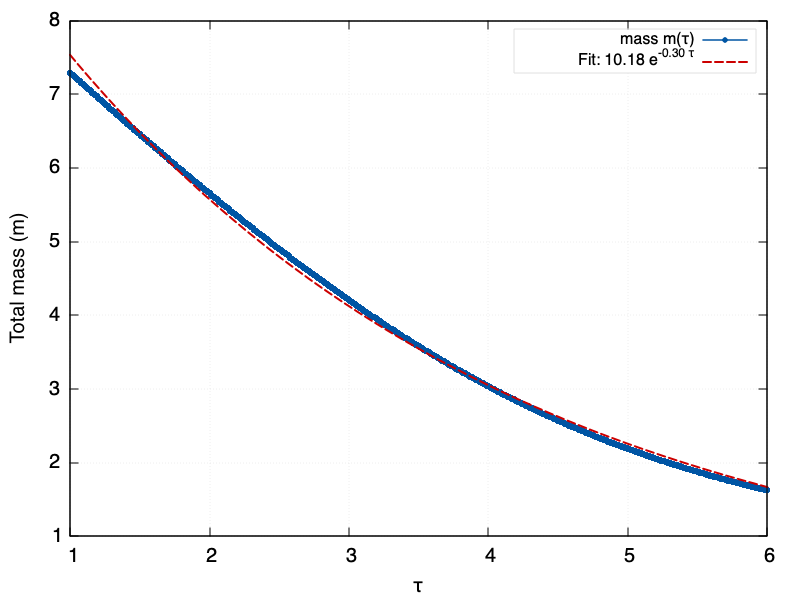}
        \caption{Evolution of the mass of the stationary branch as a function of $\tau$ when $\mu=1$ - Exponential decrease}
        \label{fig:mass_exponential}
    \end{subfigure}\hfill
    \begin{subfigure}[t]{0.48\textwidth}
        \centering
        \includegraphics[width=\linewidth]{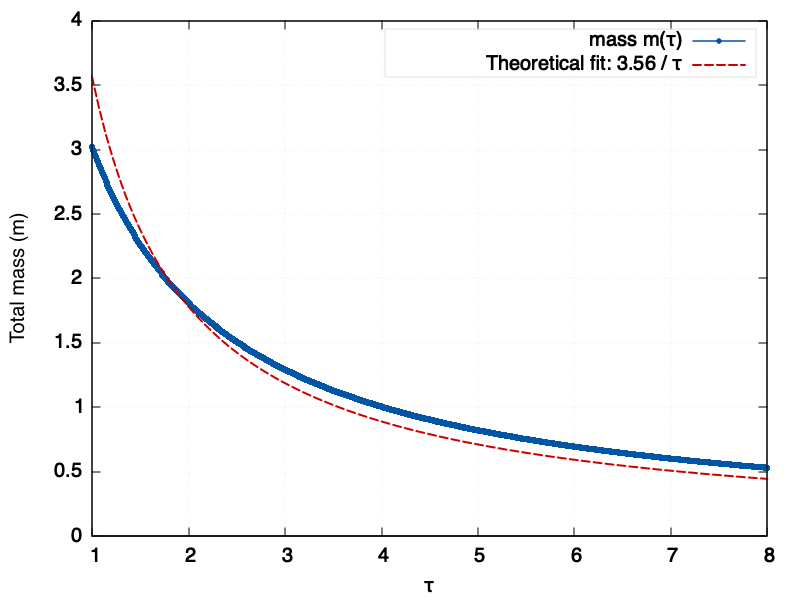}
        \caption{Evolution of the mass of the stationary branch as a function of $\tau$ when $\mu=0$ - Hyperbolic decrease}
        \label{fig:mass_hyperbolic}
    \end{subfigure}
    \vspace{-2mm}
    \caption{Evolution of the mass of the stationary branch as a function of $\tau$ when $\beta=1$.}
    \label{fig:mass_evolution_beta>0}
\end{figure}

As illustrated in Figure \ref{fig:stability_analysis_beta>0}, the real parts of all eigenvalues remain strictly negative for both $\mu=1$ and $\mu=0$, confirming that the positive stationary state is linearly stable without condition when $\beta > 0$. Furthermore, tracking the spectrum as $\mu \to 0$ at fixed parameters $\tau$ and $\beta$ (Figure \ref{fig:stability_analysis_mu_to_0}) shows a continuous deformation of the eigenvalues without any crossing of the imaginary axis, ensuring robustness of this stability.
Figure \ref{fig:mass_evolution_beta>0} highlights the evolution of the total mass $m(\tau)$ of the stationary branch. To accurately track this behaviour, for each parameter value $\tau$ along the computed branch, we extract the exact stationary profile from the Newton--Raphson solver and compute its corresponding total mass $m(\tau)$. For $\mu=1$, the mass undergoes an exponential decay (Figure \ref{fig:mass_exponential}), whereas for $\mu=0$, it follows a hyperbolic decay (Figure \ref{fig:mass_hyperbolic}). In both scenarios, the mass stays strictly bounded away from zero. The fact that this strictly positive branch extends indefinitely for all $\tau$ corroborates our previous existence result: as long as $\beta$ is strictly positive, we have the existence of a positive solution to the equation with transfer. So, regardless of the values of $\tau$ and $\mu$, the mass will always be positive. This is also consistent with the stability analysis of the zero-population state: with an infinitesimal population, both transfer and competition are neglected, and so there is no extinction.
\newpage
\subsubsection{When $\beta=0$}
We set $\mu = 1$ and recall that a sufficient analytical condition for the existence of a positive stationary solution is:

$$\frac{\tau}{\mu} < \lambda_0 \quad (\text{where } \lambda_0 \text{ is the principal eigenvalue of the transfer-free equation}).$$

Better estimates for \eqref{eq:debut_priori} could provide more accurate conditions. 
Numerically, we see that there exist a $\lambda_c\geq\lambda_0$ such that if:

$$\frac{\tau}{\mu} < \lambda_c,$$

a stationary state exists and is stable, while the zero state is unstable.
\begin{figure}[h!]
    \centering
    \begin{subfigure}[t]{0.48\textwidth}
        \centering
        \includegraphics[width=\linewidth]{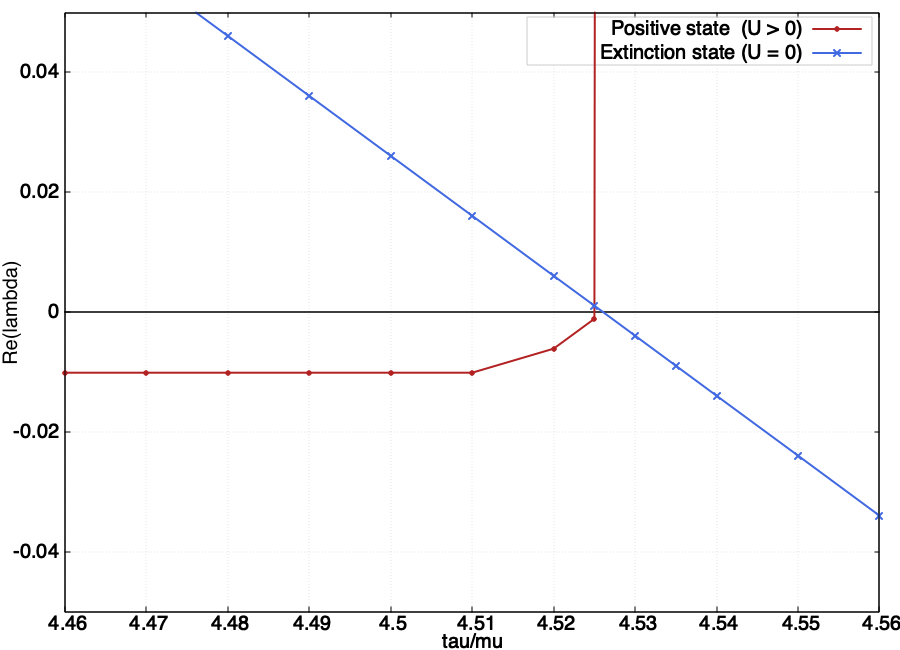}
        \caption{Real part of the leading eigenvalue for the stationary and zero states as a function of $\tau/\mu$. The exchange of stability happens for $\lambda_c \simeq 4.525$. Here $\lambda_0 \simeq 4.48$.}
        \label{fig:dominant_eigenvalue_beta_0}
    \end{subfigure}
    \hfill
    \begin{subfigure}[t]{0.48\textwidth}
        \centering
        \includegraphics[width=\linewidth]{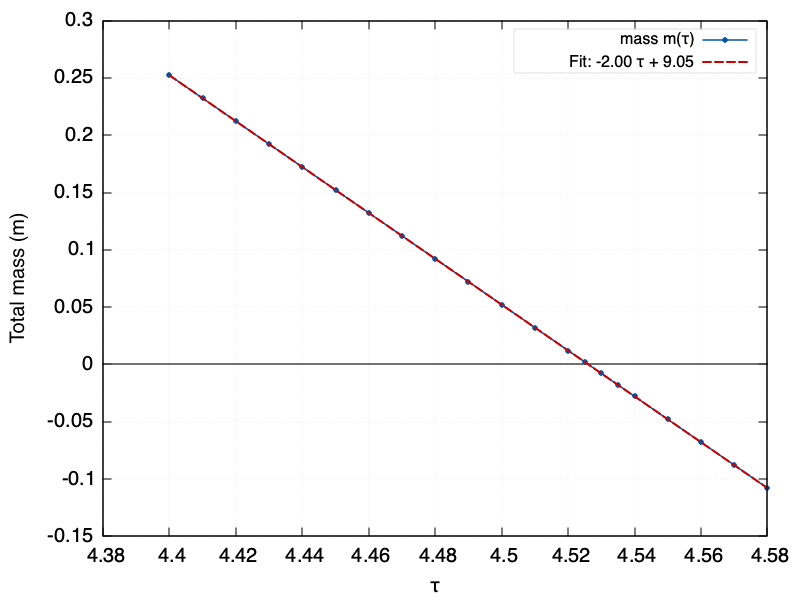}
        \caption{Evolution of the mass of the stationary solution as a function of $\tau$ when $\beta=0$, $\mu=1$.}
        \label{fig:mass_beta_0_fd}
    \end{subfigure}

    \vspace{0.8em} 

    \begin{subfigure}[t]{0.48\textwidth}
        \centering
        \includegraphics[width=\linewidth]{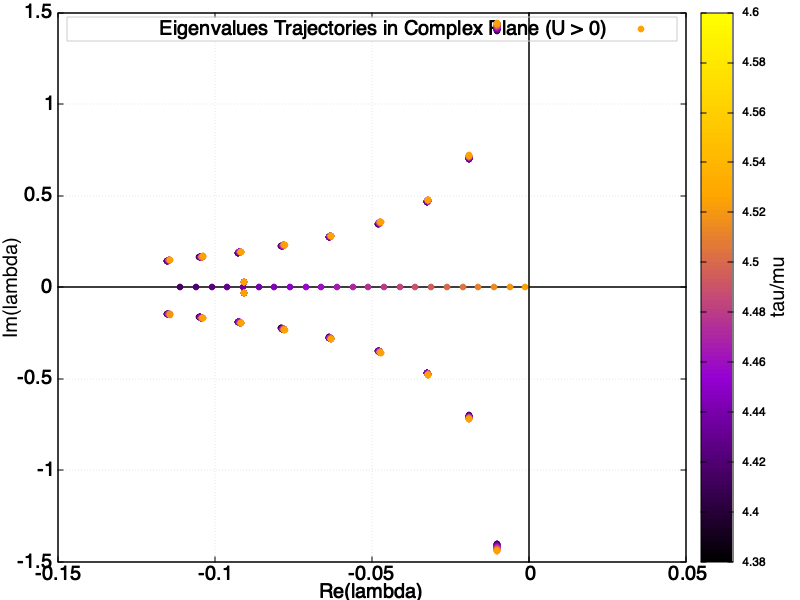}
        \caption{Spectrum evolution of the Jacobian as a function of $\tau/\mu$ - $\beta=0$ - Zoom on the left quadrant.}
        \label{fig:jacobian_spectrum_beta_0}
    \end{subfigure}
    \hfill
    \begin{subfigure}[t]{0.48\textwidth}
        \centering
        \includegraphics[width=\linewidth]{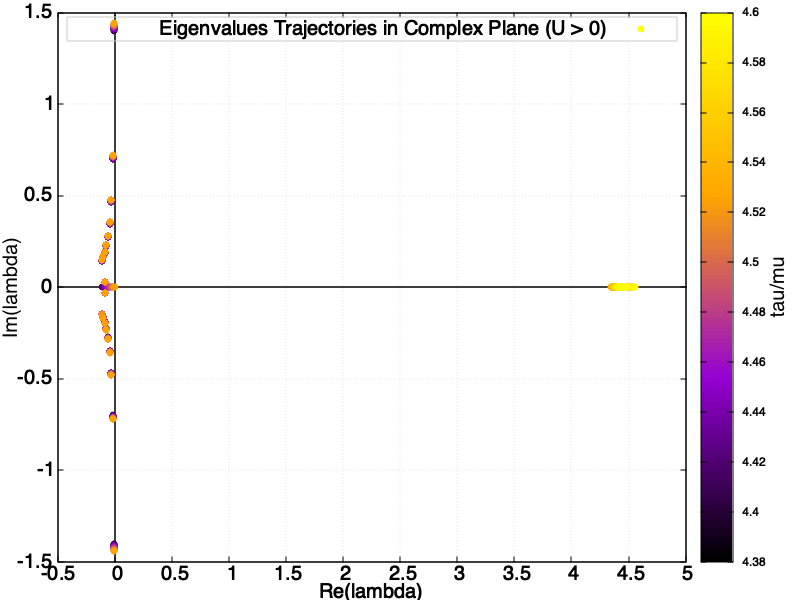}
        \caption{Spectrum of the Jacobian as a function of $\tau/\mu$ - $\beta=0$ - Full picture including the leading real eigenvalues that cross the imaginary axis.}
        \label{fig:mass_beta_0_partie_plus}
    \end{subfigure}

    \caption{Global analysis of eigenvalues, spectra, and stationary mass evolution as a function of $\tau$ and $\tau/\mu$.}
    \label{fig:combined_stability_analysis}
\end{figure}
\newpage
Once $\frac{\tau}{\mu} > \lambda_c$, the positive stationary state ceases to exist and the zero state becomes stable (in practice, the stationary solutions persists but corresponds to a state with negative mass).
Figure \ref{fig:mass_beta_0_fd} demonstrates that the mass of the stationary branch decreases linearly as a function of $\tau$ and crosses into negative, non-physical values at the critical threshold $\tau/\mu \approx 4.525$. This geometric crossing coincides exactly with the spectral exchange of stability observed in Figure \ref{fig:dominant_eigenvalue_beta_0}, where the leading eigenvalue of the stationary state becomes positive while the zero state's leading eigenvalue shifts to negative.

As $\tau$ increases, the branch of solutions crosses the horizontal axis. This transition corresponds to a transcritical bifurcation (see \cite{kuznetsov1998elements}) at $\tau/\mu \approx \lambda_c $, where the positive branch crosses the trivial state and exchanges its stability. Since the critical eigenvalue is purely real ($\operatorname{Im}(\lambda) = 0$), no local Hopf bifurcation takes place, ruling out the emergence of small-amplitude limit cycles near this threshold. To be synthetic, we have

\begin{center}
\renewcommand{\arraystretch}{1.3}
\begin{tabular}{|l|r@{\extracolsep{0pt}}|l|}
\hline
$\tau$ & \multicolumn{1}{l}{$0$} \hfill & $\lambda_c$ \hfill $+\infty\ $ \\
\hline
\textbf{Stationary state $> 0$} & \multicolumn{1}{c|}{stable} & \multicolumn{1}{c|}{ceases to exist} \\
\hline
\textbf{Zero state} & \multicolumn{1}{c|}{unstable} & \multicolumn{1}{c|}{stable} \\
\hline
\end{tabular}
\end{center}

\vspace{1em}

\noindent\textbf{Remark: Transient oscillations and absence of limit cycles.} 
Initial observations of the integro-differential equation in the frequency-dependent regime ($\beta=0$) suggested the potential existence of undamped limit cycles for certain parameter ranges. However, extended time integration reveals this to be a misconception caused by extremely slow relaxation dynamics.

\begin{figure}[h!]
    \centering
    \begin{subfigure}[t]{0.45\textwidth}
        \centering
        \includegraphics[width=\linewidth]{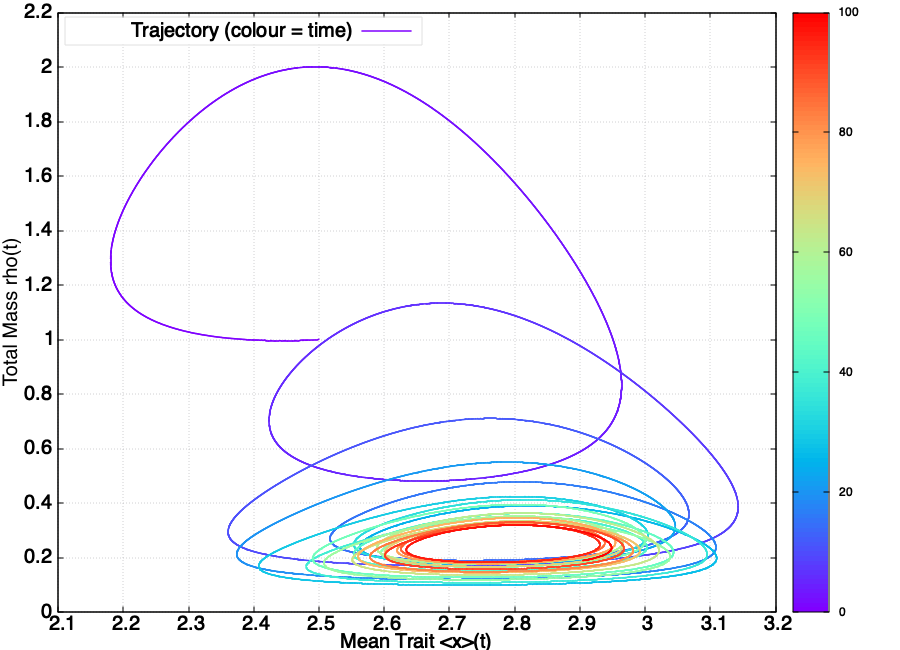}
        \caption{Time $t \in [0, 100]$}
    \end{subfigure}
    \hfill
    \begin{subfigure}[t]{0.48\textwidth}
        \centering
        \includegraphics[width=\linewidth]{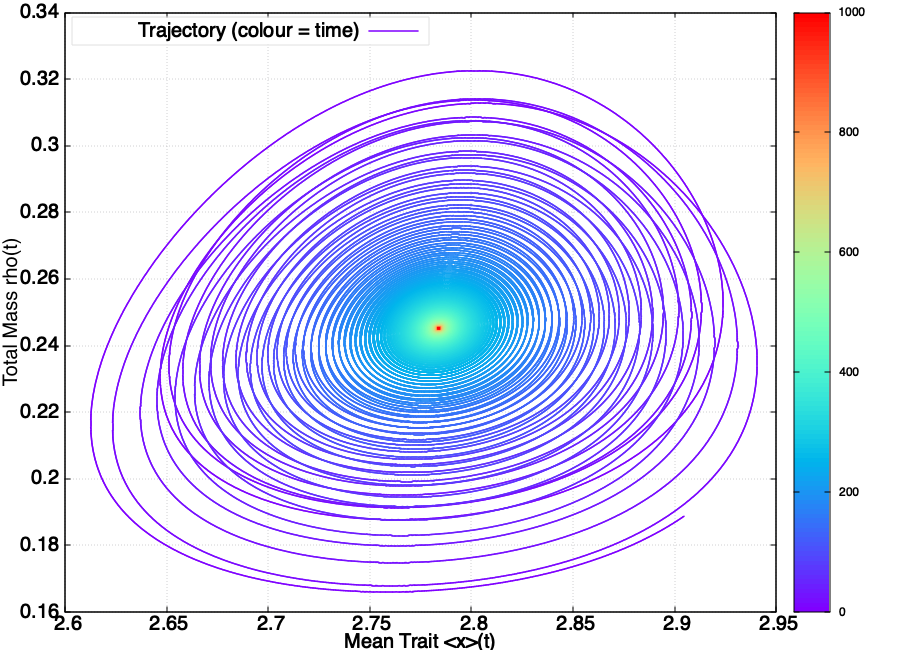}
        \caption{Time $t \in [100, 1100]$}
    \end{subfigure}
    \caption{Phase space orbit (Total mass vs. Mean trait) from a uniform initial condition, showing the slow convergence to a stable focus ($\beta=0, \tau=4.4$) - The time intervals $[0,100]$ and $[100,1100]$ are represented separately..}
    \label{fig:orbites_uniform}
\end{figure}

\begin{figure}[h!]
    \centering
    \begin{subfigure}[t]{0.48\textwidth}
        \centering
        \includegraphics[width=\linewidth]{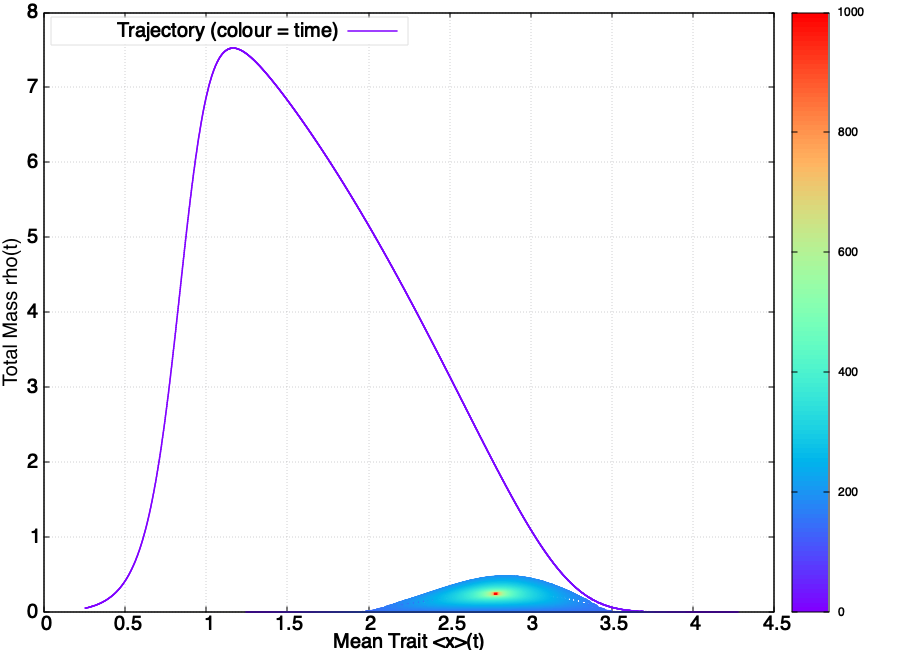}
        \caption{Full trajectory}
    \end{subfigure}
    \hfill
    \begin{subfigure}[t]{0.48\textwidth}
        \centering
        \includegraphics[width=\linewidth]{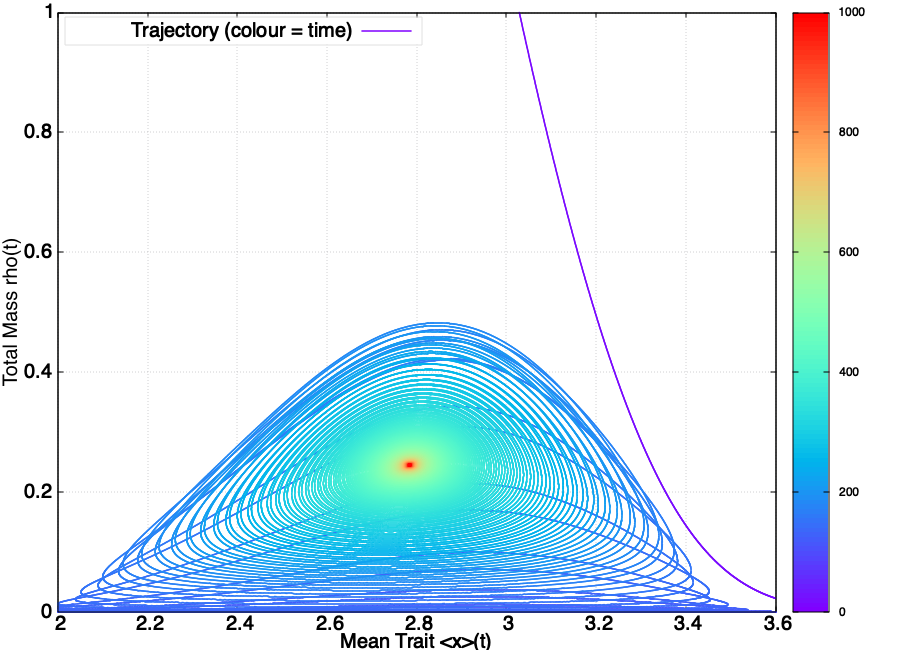}
        \caption{Zoom on the stable focus}
    \end{subfigure}
    \caption{Phase space orbit from a Gaussian initial condition centered at 0 ($\beta=0, \tau=4.4$)}
    \label{fig:orbites_gaussian}
\end{figure}

To correctly visualize this asymptotic behaviour, we simulated the system using the forward Euler scheme and projected the solution onto a 2D phase space, plotting the total population mass against the mean trait. As shown in Figures \ref{fig:orbites_uniform} and \ref{fig:orbites_gaussian}, starting from either a uniform or a Gaussian initial condition, the trajectories exhibit slow, inward-spiraling oscillations. The system ultimately settles at a stable focus, requiring considerable simulation time (up to $t=1100$) to fully dampen the transients. This confirms that as long as $\tau/\mu < \lambda_c$, the population robustly converges to the stationary equilibrium, ruling out the existence of permanent limit cycles even at the singular limit $\beta=0$.

To conclude, our numerical investigation highlights a contrast in the population dynamics depending on the transfer regime. In the regular case where $\beta > 0$, the system robustly and exponentially converges toward a positive stationary state. This non-trivial profile remains unconditionally linearly stable, which perfectly complements the analytical instability of the zero state. Conversely, the singular limit $\beta = 0$ fundamentally alters this behaviour. The system undergoes a transcritical bifurcation at the critical threshold $\tau/\mu \approx \lambda_c$, where the positive branch crosses the trivial state and exchanges its stability. Finally, despite the emergence of long-lived oscillations near this threshold, extended phase-space tracking confirms that these are purely transient structures mapping to a stable focus. The population systematically relaxes to a fixed equilibrium, definitively ruling out the existence of permanent limit cycles within the explored parameter space.

\section*{Acknowledgments}
This work has been supported by the Chair “Modélisation Mathématique
et Biodiversité” of Veolia Environnement-Ecole Polytechnique-Museum National d’Histoire
Naturelle-Fondation X. Funded by the European Union (ERC, SINGER, 101054787). Views and opinions expressed are those of the author(s) only and do not necessarily reflect those of the European Union or the European Research Council. Neither the European Union nor the granting authority can be held responsible for them. The author would like to thank Sylvie Méléard and Viet Chi Tran for their supervision, and Pierre Collet for his kind analytical and numerical insights.

\vspace{3mm}
\noindent{\bf AI Statement:} Gemini (Google) was utilized solely for numerical code optimization, LaTeX typesetting and layout assistance, and language translation. All mathematical proofs and core theoretical developments were devised entirely by the author without AI tools. The author takes full responsibility for the final manuscript.
\nocite{*}
\bibliography{biblio}
\begin{appendix}
\section{Useful results and definitions}\label{appendix}
\begin{theorem}[Dini]\label{thm:dini}
Let $(f_n)$ be a sequence of functions defined on the interval $[a,b]$ and taking values in $\mathbb{R}$, and let $f:[a,b]\mapsto \mathbb{R}$ be continuous. If each function $f_n$ is increasing and $f_n$ converges simply to $f$ on $[a,b]$, then the convergence is uniform.
\end{theorem}
For the compactness of operators in the space $L^1(d\nu)$, we make use of two classical Theorems which can be found in Chapter IV of \cite{brezis2011functional}.
We denote the translations
\begin{equation*}
    \tau_h f(x)=f(x+h), x\in\mathbb{R}^N,h \in\mathbb{R}^N.
\end{equation*}
\begin{theorem}[Kolmogorov-Riesz-Fréchet]\label{thm:kolmogorov}

Let $\mathcal{F}$ be a bounded set in $L^p(\mathbb{R}^N)$ with $1\leq p<\infty$. If 
\begin{equation}\label{eq:continuite_translations}
    \forall \varepsilon>0,\,\exists \delta>0,\,\|\tau_h f-f\|_p<\varepsilon,\, \forall f \in \mathcal{F}, \forall h\in\mathbb{R}^N\, \text{with }|h|<\delta
\end{equation}
and that for every $\varepsilon>0$, there exists a bounded measurable subset $\Omega\subset\mathbb{R}^N$ such that 
\begin{equation*}
    \|f\|_{L^p(\mathbb{R}^N\backslash \Omega)}<\varepsilon, \forall f\in\mathcal{F}
\end{equation*}
then $\mathcal{F}$ is relatively compact in $L^p(\mathbb{R}^N)$.
\end{theorem}
\begin{theorem}[Arzela-Ascoli]\label{thm:arzela_ascoli}

Let $K$ be a compact metric space and $\mathcal{H}$ a bounded subset of $\mathcal{C}(K)$. We assume that $\mathcal{H}$ is uniformly equicontinuous, that is to say
\begin{equation*}
    \forall\varepsilon,\,\exists\delta>0,\, d(x_1,x_2)<\delta \implies |f(x_1)-f(x_2)|<\varepsilon,\, \forall f\in\mathcal{H}.
\end{equation*}
Then $\mathcal{H}$ is relatively compact in $\mathcal{C}(K)$.
\end {theorem}
For the solution to the eigenvalue problem, we use Jentzsch’s theorem, which can be found in Part V, Chapter 6 of \cite{schaefer1974banach}.
\begin{theorem}[Jentzsch]\label{thm:Jentzsch}

Let $E:=L^p(d\mu)$, with $1\leq p\leq \infty$ and $(X,\Sigma,\mu)$ a $\sigma$-finite measured space. Let us consider $T$ a linear operator given by a kernel $K\geq0$, which is $(\Sigma\times\Sigma)$-measurable :
\begin{equation*}
    \forall f\in E,\, Tf(t):=\int_X K(s,t)f(s)\mu(ds).
\end{equation*}

Suppose that
\begin{itemize}
    \item A certain power of $T$ is compact.
    \item For $S\in\Sigma$ such that $\mu(S)>0,\,\mu(X\backslash S)>0$, we have 
    \begin{equation*}\label{Jentzsch_s_property}
        \int_{X\backslash S}\int_S K(s,t)\mu(ds)\mu(dt)>0.
    \end{equation*}
    Then $s(T)$ is an eigenvalue of $T$ with a unique renormalised eigenfunction $f$ satisfying $f(s)>0$ for $\mu$-almost all $s$. Furthermore, if $K(s,t)>0$ $(\mu\otimes\mu)$-almost everywhere, then every other eigenvalue $\lambda$ of $T$ has a modulus such that $|\lambda|<s(T)$.
\end{itemize}
\end{theorem}

For the continuity of the spectral radius, we use a result from \cite{kato2013perturbation}, Section 3.5 of Chapter IV.
\begin{theorem}[Continuity of the spectral radius]\label{thm:continuite_du_rayon_spectral}

Let $I$ be an open interval of $\mathbb{R}$ and $\lambda$ an element of $I$. Let $T_\lambda$ denote a compact operator on $L^1(d\nu)$.
If $\lambda\mapsto T_\lambda$ is continuous in the space of bounded operators, that is to say
\begin{equation*}
    \lim_{\lambda\rightarrow\lambda_0} \sup_{\|f\|_\nu\leq 1}\|T_\lambda f-T_{\lambda_0}f\|_\nu=0
\end{equation*}
then every simple eigenvalue $\psi(\lambda_0)$ of $T_{\lambda_0}$ can be extended to a continuous function $\lambda \mapsto \psi(\lambda)$ on a neighbourhood of $\lambda_0$. Furthermore, if the eigenvalue is simple for every admissible $\lambda$, then $\psi$ is as regular as $\lambda\longmapsto T_\lambda$. 
Under the same conditions and assuming the invertibility of $\lambda\rightarrow\psi(\lambda)$, the mapping which associates the corresponding normalised eigenvector with $\psi(\lambda)$ is continuous.
\end{theorem}
For the strict decrease of the spectral radius, we use the generalisation of the Perron-Frobenius theorem and the Collatz-Wielandt formula to positive operators on Banach spaces, established by Marek \cite{Marek1970frobenius}.

\begin{definition}[Semi-nonsupporting operator]\label{def:smnnsp}
Let $E$ be a Banach space ordered by a closed convex cone $K$. A bounded, positive linear operator $T: E \to E$ is said to be \textit{semi-nonsupporting} if, for every $u \in K \setminus \{0\}$ and for every non-zero positive linear form $u' \in K' \setminus \{0\}$, there exists an integer $p \ge 1$ such that:
\begin{equation*}
    \langle T^p u, u' \rangle > 0.
\end{equation*}
\end{definition}
\begin{definition}[Reproducing cone]

A cone $K$ is said to be reproducing if and only if every element of the space $E$ can be written as the difference of two elements of the cone:
\begin{equation*}
    E=K-K
\end{equation*}
\end{definition}
\begin{definition}[Normal cone]
A cone $K$ is said to be normal if there exists a constant $\delta>0$ such that for all $x,y\in K$:
\begin{equation*}
    \|x+y\|\geq \delta\|x\|.
\end{equation*}
\end{definition}
We state a Lemma regarding the existence of the dual eigenvector, which is necessary for the application of Marek’s $Min-Max$ formula.
\begin{lemma}[Existence of the dual eigenvector \cite{Marek1970frobenius}, Thm 2.2]\label{lemma:marek}
Let $K$ be a cone in $E$ and $T$ be a linear operator. If $K$ is stable under $T$, then the spectral radius $s(T)$ is an element of the spectrum of $T$. 
If $T$ is compact, the spectral radius is a pole of the resolvent. There exists at least one $g>0$ belonging to $K'$, a cone of the dual space $E'$, such that $T' g=s(T)g$. 
\end{lemma}
\begin{theorem}[Marek’s Min-Max Formula \cite{Marek1970frobenius}, Thm 3.2]\label{thm:Marek}
Let $E$ be a Banach space equipped with a closed, normal and reproducing cone $K$. Let $T$ be a bounded, positive and semi-nonsupporting linear operator.
Suppose further that the spectral radius $s(T)$ is an eigenvalue of $T$ associated with an eigenvector in $K$, and an eigenvalue of the adjoint $T'$ associated with an eigenform in $K'$.

Then, the spectral radius $s(T)$ is characterised by the following formula:
\begin{equation*}\label{eq:collatz_wielandt}
    s(T) = \sup_{u \in K \setminus \{0\}} r_u(T) = \inf_{u \in K \setminus \{0\}} r^u(T),
\end{equation*}
where the local lower and upper bounds are defined by:
\begin{align*}
    r_u(T) &= \sup \{ \nu \in \mathbb{R} : Tu \ge \nu u \}, \\
    r^u(T) &= \inf \{ \tau \in \mathbb{R} : Tu \le \tau u \}.
\end{align*}
\end{theorem}

To obtain a sufficient condition and a biological intuition regarding the parameters underlying the existence of the no-transfer problem, in other words, to estimate $\lambda_0$, we require a result concerning the largest eigenvalue of an upper-bounded self-adjoint operator. Such a result can be found in Theorem XIII.1 of \cite{reed1978iv}
\begin{theorem}[Min-Max Principle – Variational Formula]\label{thm:min_max_variationnel}.

Let $H$ be a self-adjoint operator in a space $E$ and lower-bounded $( H\geq cI)$. We define
\begin{equation*}
    \mu_1(H)=\inf_{\phi\in D(H), \|\phi\|} \langle H\phi,\phi\rangle.
\end{equation*}
This quantity corresponds to one of the following two situations:
\begin{itemize}
    \item If $\mu_1$ lies below the essential spectrum, then it is the smallest eigenvalue of the operator.
    \item If $\mu_1$ is not an isolated eigenvalue at the bottom of the spectrum, then $(\mu_1=\inf \sigma_{ess}(H))$.
\end{itemize}
\end{theorem}
If we have an operator that is bounded above and self-adjoint, we take the $\sup$. Furthermore, if we know that the largest eigenvalue of our operator is isolated, then its value is given by $\mu_1$.
Finally, the fixed-point theorem we use is the Leray–Schauder theorem \cite{leray1934topologie, mawhin1999leray}. 
\begin{theorem}[Leray–Schauder]\label{thm:leray_schauder}

Let $X$ be a Banach space and $I=[0,1]$. Let $\mathcal{A}\subset X$ be a bounded open set with closure $\overline{\mathcal{A}}$ and boundary $\partial A$, and let $\Gamma:\overline{\mathcal{A}}\times I\rightarrow X$ be an operator. Let $F_\Gamma$ denote the set of fixed points (which may be empty) defined by 
\begin{equation*}
    F_\Gamma=\{(x,\delta)\in\overline{\mathcal{A}}\times I|\,x=\Gamma(x,\delta)\}.
\end{equation*}
We assume:
\begin{itemize}
    \item $\Gamma:\overline{\mathcal{A}}\times I\rightarrow X$ is continuous and compact.
    \item $F_\Gamma\cap \partial\mathcal{A}=\emptyset$, in other words, no fixed point in $\overline{\mathcal{A}}$ belongs to the boundary of $\overline{\mathcal{A}}$.
    \item The operator $\Gamma(\cdot,0)$ has a unique non-zero fixed point.
\end{itemize}
Then, for every $\delta$ in $[0,1]$, there exists an $x$ in $X$ such that $\Gamma(x,\delta)=x$. In particular, the equation $\Gamma(x,1)=x$ has at least one solution.
\end{theorem}

\end{appendix}
\end{document}